\documentclass[12pt]{amsart}
\usepackage{amssymb}
\usepackage{pstricks,pstricks-add}

\newtheorem{theorem}{Theorem}[section]

\newtheorem{proposition}[theorem]{Proposition}
\newtheorem{corollary}[theorem]{Corollary}
\newtheorem{definition}[theorem]{Definition}
\newtheorem{conj}[theorem]{Conjecture}

\theoremstyle{definition}
\newtheorem{remark}[theorem]{Remark}
\newtheorem{example}[theorem]{Example}

\numberwithin{equation}{section}

\newcommand{\Z}{\mathbb{Z}}
\newcommand{\R}{\mathbb{R}}
\newcommand{\C}{\mathbb{C}}
\newcommand{\Q}{\mathbb{Q}}
\newcommand{\K}{\mathbb{K}}

\newcommand{\F}{\mathcal{F}}

\newcommand{\dbar}{\bar\partial}

\title{A beginner's introduction to Fukaya categories}
\author{Denis Auroux}
\address{Department of Mathematics, UC Berkeley, Berkeley CA 94720-3840, USA}
\email{auroux@math.berkeley.edu}
\thanks{The author was partially supported by NSF grant DMS-1007177}

\begin{document}
\maketitle

This text is based on a series of lectures given at a Summer School on
Contact and Symplectic Topology at Universit\'e de Nantes in June 2011. 

The goal of these notes is to give a short introduction to Fukaya
categories and some of their applications. The first half of the text
is devoted to a brief review of Lagrangian Floer (co)homology and product
structures.  Then we introduce the Fukaya category (informally and 
without a lot of the necessary technical detail), and briefly discuss 
algebraic concepts such as exact triangles and generators. Finally, 
we mention wrapped Fukaya categories and 
outline a few applications to symplectic 
topology, mirror symmetry and low-dimensional topology.

These notes are in no way a comprehensive text on the subject; however
we hope that they will provide a useful introduction to Paul Seidel's book
\cite{SeBook} and other texts on Floer homology, Fukaya categories, and
their applications. We assume that the reader is generally familiar with 
the basics of symplectic geometry, and some prior exposure to
pseudo-holomorphic curves is also helpful; the reader is
referred to \cite{McS1,McS2} for background material.

\subsection*{Acknowledgements}
The author wishes to thank the organizers of the Nantes Trimester
on Contact and Symplectic Topology for the pleasant atmosphere
at the Summer School, and Ailsa Keating for providing a copy of the
excellent notes she took during the lectures. Much of the material
presented here I initially learned from Paul Seidel and Mohammed 
Abouzaid, whom I thank for their superbly written papers and their
patient explanations. Finally, the author was partially supported
by an NSF grant (DMS-1007177).

\section{Lagrangian Floer (co)homology}\label{s:HF}

\subsection{Motivation}
Lagrangian Floer homology was introduced by Floer in the late 1980s in order
to study the intersection properties of compact Lagrangian submanifolds in
symplectic manifolds and prove an important case of Arnold's conjecture
concerning intersections between Hamiltonian isotopic Lagrangian
submanifolds \cite{Arnold}.

Specifically, let $(M,\omega)$ be a symplectic manifold (compact, or
satisfying a ``bounded geometry'' assumption), and let $L$ be a compact Lagrangian submanifold of $M$.
Let $\psi\in \mathrm{Ham}(M,\omega)$ be a Hamiltonian diffeomorphism.
(Recall that a time-dependent Hamiltonian $H\in C^\infty(M\times
[0,1],\R)$ determines a family of
Hamiltonian vector fields $X_t$ via the equation
$\omega(\cdot,X_t)=dH_t$, where $H_t=H(\cdot,t)$;
integrating these vector fields over $t\in [0,1]$ yields the Hamiltonian
diffeomorphism $\psi$ {\em generated} by $H$.)  

\begin{theorem}[Floer \cite{Floer}]\label{thm:floerarnold}
Assume that the symplectic area of any topological disc in $M$ with boundary 
in $L$ vanishes. Assume moreover that $\psi(L)$ and $L$ intersect
transversely. Then the number of intersection points of $L$ and $\psi(L)$
satisfies the lower bound $|\psi(L)\cap L|\ge \sum_i \dim H^i(L;\Z_2)$.
\end{theorem}

\noindent
Note that, by Stokes' theorem, since $\omega_{|L}=0$, the symplectic area
of a disc with boundary on $L$ only depends on its class in the relative
homotopy group $\pi_2(M,L)$.

The bound given by Theorem \ref{thm:floerarnold} is stronger than what one could
expect from purely topological considerations. The assumptions that
the diffeomorphism $\psi$ is Hamiltonian, and that $L$ does not bound discs of
positive symplectic area, are both essential (though the latter can be
slightly relaxed in various ways).

\begin{figure}[t]
\setlength{\unitlength}{5mm}
\begin{picture}(4,6)(-2,-3)
\newgray{gray30}{0.8}
\newgray{gray15}{0.9}
\psset{unit=\unitlength}
\psellipticarc(0,-2.5)(2,0.5){-180}{0}
\psellipticarc[linestyle=dashed](0,-2.5)(2,0.5){0}{180}
\psellipse(0,2.5)(2,0.5)
\psline(-2,-2.5)(-2,2.5)
\psline(2,-2.5)(2,2.5)
\psccurve[linestyle=dashed,fillstyle=solid,fillcolor=gray15]%
  (-2,0)(-2,0)(-2,0.2)(-1,0.6)(1,0.6)(2,0.2)(2,0)(2,0)%
  (2,-1)(2,-1)(2,-0.8)(1,0.1)(-1,1.2)(-2,1.2)(-2,1)(-2,1)
\psccurve[fillstyle=solid,fillcolor=gray30]%
  (-2,0.1)(-2,0.1)(-2,-0.1)(-1,-0.5)(1,-0.5)(2,-0.1)(2,0.1)(2,0.1)%
  (2,-0.9)(2,-0.9)(2,-1.1)(1,-1.1)(-1,0)(-2,0.9)(-2,1.1)(-2,1.1)
\psccurve[linestyle=dashed]%
  (-2,0)(-2,0)(-2,0.2)(-1,0.6)(1,0.6)(2,0.2)(2,0)(2,0)%
  (2,-1)(2,-1)(2,-0.8)(1,0.1)(-1,1.2)(-2,1.2)(-2,1)(-2,1)
\pscircle*(-0.15,-0.6){0.12}
\pscircle*(0.15,0.7){0.12}
\put(-2.4,0){\makebox(0,0)[rc]{$L$}}
\put(-2.2,1.2){\makebox(0,0)[rc]{$\psi(L)$}}
\put(-0.5,-1.1){\small $p$}
\put(0.3,1.1){\small $q$}
\end{picture}
\qquad\qquad
\begin{picture}(4,6)(-2,-3)
\newgray{gray30}{0.8}
\newgray{gray15}{0.9}
\psset{unit=\unitlength}
\psellipticarc(0,-2.5)(2,0.5){-180}{0}
\psellipticarc[linestyle=dashed](0,-2.5)(2,0.5){0}{180}
\psellipse(0,2.5)(2,0.5)
\psline(-2,-2.5)(-2,2.5)
\psline(2,-2.5)(2,2.5)
\psccurve[fillstyle=solid,fillcolor=gray30]
  (-1.5,0.5)(0,1.5)(1,1)(0,0.5)(-1,-0.5)
\psccurve[fillstyle=solid,fillcolor=gray15]
  (1.5,-1.2)(0,-2.2)(-1,-1.7)(0,-1.2)(1,-0.2)
\put(-1.3,1.4){\makebox(0,0)[rc]{$L$}}
\put(1.5,-0.5){\makebox(0,0)[lc]{$\psi(L)$}}
\end{picture}
\caption{Arnold's conjecture on the cylinder $\R\times S^1$: an example
(left) and a non-example (right)}
\label{fig:T*S1}
\end{figure}
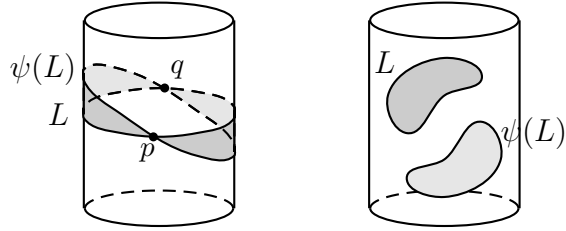

\begin{example}
Consider the cylinder $M=\R\times S^1$, with the standard area form, and a
simple closed curve $L$ that goes around the cylinder once: then $\psi(L)$
is also a simple closed curve going around the cylinder once, and the 
assumption that $\psi\in \mathrm{Ham}(M)$ means that the total signed area
of the 2-chain bounded by $L$ and $\psi(L)$ is zero.  It is then an
elementary fact that $|\psi(L)\cap L|\ge 2$, as claimed by Theorem
\ref{thm:floerarnold}; see Figure~\ref{fig:T*S1} left. On the other hand, the
result becomes false if we only assume that $\psi$ is a symplectomorphism
(a large vertical translation of the cylinder is area-preserving and
eventually displaces $L$ away from itself); or if we take $L$ to be a
homotopically trivial simple closed curve, which bounds a disc of
positive area (see Figure~\ref{fig:T*S1} right).
\end{example}

Floer's approach is to associate to the pair of Lagrangians
$(L_0,L_1)=(L,\psi(L))$ a
chain complex $CF(L_0,L_1)$,
freely generated by the intersection points of $L_0$ and $L_1$, equipped
with a differential $\partial:CF(L_0,L_1)\to CF(L_0,L_1)$, with
the following properties:
\begin{enumerate}
\item $\partial^2=0$, so the {\em Floer cohomology}
$HF(L_0,L_1)=\mathrm{Ker}\,\partial/\mathrm{Im}\,\partial$ is
well-defined;
\item if $L_1$ and $L'_1$ are Hamiltonian isotopic then $HF(L_0,L_1)\simeq
HF(L_0,L'_1)$;
\item if $L_1$ is Hamiltonian isotopic to $L_0$, then
$HF(L_0,L_1)\simeq H^*(L_0)$ (with suitable coefficients).
\end{enumerate}
Theorem \ref{thm:floerarnold} then follows immediately, since the rank of
$HF(L,\psi(L))\simeq H^*(L)$ is bounded by that
of the Floer complex $CF(L,\psi(L))$, which equals $|\psi(L)\cap L|$.

Formally, Lagrangian Floer (co)homology can be viewed as an infinite-dimensional 
analogue of Morse (co)homology for the {\em action functional} on (the universal
cover of) the path space $\mathcal{P}(L_0,L_1)=\{\gamma:[0,1]\to M\,|\,\gamma(0)\in L_0,\ 
\gamma(1)\in L_1)\}$, 
$$\mathcal{A}(\gamma,[\Gamma])=-\textstyle\int_\Gamma \omega,$$
where $(\gamma,[\Gamma])\in \tilde{\mathcal{P}}(L_0,L_1)$ consists of a path
$\gamma\in \mathcal{P}(L_0,L_1)$ and an equivalence class $[\Gamma]$ of a
homotopy $\Gamma:[0,1]\times [0,1]\to M$ between $\gamma$ and a fixed base 
point in the connected component of $\mathcal{P}(L_0,L_1)$
containing $\gamma$. The critical points of $\mathcal{A}$ are (lifts of)
constant paths at intersection points, and its gradient flow lines 
(with respect to the natural $L^2$-metric induced by $\omega$ and a compatible
almost-complex structure) are pseudo-holomorphic strips bounded by $L_0$ and
$L_1$. 

However, the analytic difficulties posed by Morse theory in the 
infinite-dimensional setting are such that the actual definition of
Floer (co)homology does not rely on this interpretation: instead, the Floer
differential is defined in terms of moduli spaces of pseudo-holomorphic 
strips.

\subsection{The Floer differential}
Let $L_0,L_1$ be compact Lagrangian submanifolds of a symplectic manifold
$(M,\omega)$, and assume for now that $L_0$ and $L_1$ intersect
transversely, hence at a finite set of points. 

Before we introduce the
Floer complex and the Floer differential, a brief discussion of coefficients 
is in order. In general, Floer cohomology is defined with {\em Novikov 
coefficients} (over some base field $\K$, for example $\K=\Q$, or
$\K=\Z_2$).

\begin{definition}\label{def:novikov}
The {\em Novikov ring} over a base field $\K$ is
$$\Lambda_0=\left\{\textstyle\sum\limits_{i=0}^\infty a_i
T^{\lambda_i}\,\Big|\,
a_i\in\K,\ \lambda_i\in \R_{\ge 0},\ \lim\limits_{i\to \infty}
\lambda_i=+\infty\right\}.$$
The Novikov field $\Lambda$ is the field of fractions of $\Lambda_0$,
i.e.\
$$\Lambda=\left\{\textstyle\sum\limits_{i=0}^\infty a_i
T^{\lambda_i}\,\Big|\,
a_i\in\K,\ \lambda_i\in \R,\ \lim\limits_{i\to \infty}
\lambda_i=+\infty\right\}.$$
\end{definition}

The Floer complex is then the free $\Lambda$-module generated by
intersection points: we denote by $\mathcal{X}(L_0,L_1)=L_0\cap L_1$ the
set of generators, and set
$$CF(L_0,L_1)=\bigoplus_{p\in \mathcal{X}(L_0, L_1)} \Lambda\cdot p.$$
Equip $M$ with an $\omega$-compatible almost-complex structure
$J$. (By a classical result, the space of $\omega$-compatible
almost-complex structures
$\mathcal{J}(M,\omega)=\{J\in \mathrm{End}(TM)\,|\,J^2=-1\text{ and }
g_J=\omega(\cdot,J\cdot)\text{ is a Riemannian metric}\}$
is non-empty and contractible \cite{McS1}.)

The Floer differential $\partial:CF(L_0,L_1)\to CF(L_0,L_1)$ is defined by
counting pseudo-holomorphic strips in $M$ with boundary in $L_0$ and $L_1$:
namely, given intersection points $p,q\in \mathcal{X}(L_0, L_1)$, the coefficient of
$q$ in $\partial p$ is obtained by considering the space of
maps $u:\R\times [0,1]\to M$ which solve the Cauchy-Riemann equation
$\dbar_J u=0$, i.e.\ 
\begin{equation}\label{eq:dbar}
\dfrac{\partial u}{\partial s}+J(u)\,\dfrac{\partial u}{\partial t}=0,
\end{equation}
subject to the boundary conditions
\begin{equation}\label{eq:boundary}
\begin{cases}
u(s,0)\in L_0\text{ and } u(s,1)\in L_1 \quad \forall s\in\R,\\
\lim\limits_{s\to +\infty} u(s,t)=p,\quad \lim\limits_{s\to-\infty}
 u(s,t)=q,
\end{cases}\end{equation}
and the {\em finite energy} condition
\begin{equation}\label{eq:energy}
E(u)=\displaystyle \int u^*\omega=\iint \left|\dfrac{\partial u}{\partial s}
\right|^2\,ds\,dt<\infty.
\end{equation}

\begin{figure}[t]
\setlength{\unitlength}{8mm}%
\begin{picture}(16,3.8)(0,-1.1)
\newgray{gray30}{0.8}
\newgray{gray15}{0.9}
\psset{unit=\unitlength}
\psframe[linestyle=none,fillstyle=solid,fillcolor=gray15](-0.5,0)(4.5,1.5)
\psline[linestyle=dotted](-0.5,0)(4.5,0)
\psline[linestyle=dotted](-0.5,1.5)(4.5,1.5)
\psline(0,0)(4,0)
\psline(0,1.5)(4,1.5)
\psline{->}(0.5,0.5)(1.2,0.5)
\psline{->}(0.5,0.5)(0.5,1.2)
\put(1.1,0.65){\small $s$}
\put(0.1,1){\small $t$}
\psline{->}(5.5,0.75)(7,0.75)
\put(6.05,0.9){$u$}
\pscurve[fillstyle=solid,linestyle=none,fillcolor=gray15](8.75,0.75)(10.5,2)(12.5,2)(14.25,0.75)
\pscurve[fillstyle=solid,linestyle=none,fillcolor=gray15](8.75,0.75)(10.5,-0.5)(12.5,-0.5)(14.25,0.75)
\pscurve(7.5,-0.5)(8.75,0.75)(10.5,2)(12.5,2)(14.25,0.75)(15.5,-0.5)
\pscurve(7.5,2)(8.75,0.75)(10.5,-0.5)(12.5,-0.5)(14.25,0.75)(15.5,2)
\pscircle*(8.75,0.75){0.1}
\pscircle*(14.25,0.75){0.1}
\put(11,-1.15){$L_0$}
\put(11,2.35){$L_1$}
\put(8.6,1.15){$q$}
\put(14.2,1.15){$p$}
\end{picture}
\caption{A pseudo-holomorphic strip contributing to the Floer differential
on $CF(L_0,L_1)$}
\label{fig:strip}
\end{figure}
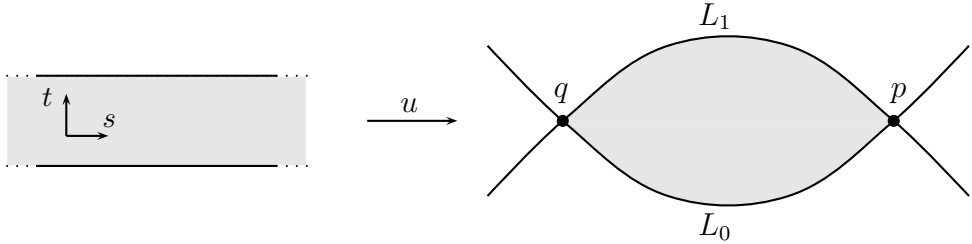

\noindent (Note that, by the Riemann mapping theorem, the strip $\R\times
[0,1]$ is biholomorphic to $D^2\setminus \{\pm 1\}$, the closed unit disc 
minus two points on its boundary; the map $u$ then extends to the closed
disc, with the boundary marked points $\pm 1$ mapping to $p$ and~$q$.)

Given a homotopy class $[u]\in \pi_2(M,L_0\cup L_1)$, we denote by 
$\widehat{\mathcal{M}}(p,q;[u],J)$ the space of solutions of
\eqref{eq:dbar}--\eqref{eq:energy} representing the class $[u]$, and
by $\mathcal{M}(p,q;[u],J)$ its quotient by the action of $\R$ by
reparametrization (i.e., $a\in \R$ acts by $u\mapsto u_a(s,t):=u(s-a,t)$).

The boundary value problem \eqref{eq:dbar}--\eqref{eq:energy} is a Fredholm
problem, i.e.\ the linearization $D_{\dbar_J,u}$ of $\dbar_J$ at a given
solution $u$ is a Fredholm operator. Specifically, $D_{\dbar_J,u}$ 
is a $\dbar$-type first-order differential operator, whose domain is a suitable
space of sections of the pullback bundle $u^*TM$ (with Lagrangian
boundary conditions), for example $W^{1,p}(\R\times [0,1],
\R\times \{0,1\};u^*TM,u_{|t=0}^*TL_0,u_{|t=1}^*TL_1)$. The Fredholm
index $\mathrm{ind}([u]):=\mathrm{ind}_\R(D_{\dbar_J,u})=\dim\mathrm{Ker}\,D_{\dbar_J,u}-
\dim\mathrm{Coker}\,D_{\dbar_J,u}$ can be computed in terms of an invariant
of the class $[u]$ called the {\em Maslov index}, which we discuss below.

The space of solutions $\widehat{\mathcal{M}}(p,q;[u],J)$ is then a smooth
manifold of dimension $\mathrm{ind}([u])$, provided that
all solutions to \eqref{eq:dbar}--\eqref{eq:energy} are
{\em regular}, i.e.\ the linearized operator $D_{\dbar_J,u}$ is surjective
at each point of $\widehat{\mathcal{M}}(p,q;[u],J)$. This transversality
property is one of three fundamental technical issues that need to be 
addressed for Floer (co)homology to be defined, the other two being the compactness 
of the moduli space $\mathcal{M}(p,q;[u],J)$, and its orientability (unless
one is content to work over $\K=\Z_2$).

Transversality and compactness will be briefly discussed in \S \ref{ss:transv}
below. On the issue of orientations, we will only consider the case where
$L_0$ and $L_1$ are oriented and spin. It is then known that
the choice of spin structures on $L_0$ and
$L_1$ determines a canonical orientation of the moduli spaces of
$J$-holomorphic strips; the construction of this orientation is fairly
technical, so we refer the reader to \cite{FO3book,SeBook} for details.

Assuming that all these issues have been taken care of, we observe that when
$\mathrm{ind}([u])=1$ the moduli space $\mathcal{M}(p,q;[u],J)$ is a compact
oriented 0-manifold, i.e.\ a finite set of points which can be counted with
signs. We can then provisionally define:

\begin{definition}\label{def:diffl}
The Floer differential $\partial:CF(L_0,L_1)\to CF(L_0,L_1)$ is the
$\Lambda$-linear map defined by 
\begin{equation}\label{eq:diffl}
\partial(p)=\sum_{\substack{q\in \mathcal{X}(L_0, L_1)\\ [u]:\,
\mathrm{ind}([u])=1}} \left(\#\mathcal{M}(p,q;[u],J)\right)\,
T^{\omega([u])}\,q,
\end{equation}
where $\#\mathcal{M}(p,q;[u],J)\in\Z$ $($or $\Z_2)$ is the signed (or
unsigned) count of points in
the moduli space of pseudo-holomorphic strips connecting $p$ to $q$ in
the class $[u]$, and $\omega([u])=\int u^*\omega$ is the symplectic
area of those strips.
\end{definition}

\noindent In general, the definition needs to be modified by introducing
a perturbation term into the Cauchy-Riemann equation in order to achieve
tranversality (see \S \ref{ss:transv} below). Thus, the Floer
differential actually 
counts {\em perturbed} pseudo-holomorphic strips connecting {\em perturbed}
intersection points of $L_0$ and $L_1$.

The following result is due to Floer for $\K=\Z_2$:

\begin{theorem}\label{thm:floer}
Assume that $[\omega]\cdot \pi_2(M,L_0)=0$ and $[\omega]\cdot
\pi_2(M,L_1)=0$. Moreover, when $\mathrm{char}(\K)\neq 2$ assume that
$L_0,L_1$ are oriented and equipped with spin structures.
Then the Floer differential $\partial$ is well-defined, satisfies $\partial^2=0$, and the Floer
cohomology
$HF(L_0,L_1)=H^*(CF(L_0,L_1),\partial)$ is, up to isomorphism, 
independent of the chosen almost-complex structure $J$ and invariant under
Hamiltonian isotopies of $L_0$ or $L_1$.
\end{theorem}

\begin{remark}
In this text we discuss the chain complex and differential for Floer {\em
cohomology}, which is dual to Floer's original construction. Namely, 
in Floer homology, the strip shown on Figure \ref{fig:strip} would 
be considered a trajectory from $q$ to $p$ rather than from $p$ to $q$,
and the grading conventions are reversed.
\end{remark}

\begin{remark} \label{rmk:novikov}
In general, the sum in the right-hand side of \eqref{eq:diffl} can be
infinite. However, Gromov's compactness theorem ensures that, given any energy
bound $E_0$, there are only finitely many homotopy classes $[u]$ with
$\omega([u])\le E_0$ for which the moduli spaces $\mathcal{M}(p,q;[u],J)$ 
are non-empty. Thus, using Novikov coefficients and weighing counts of 
strips by area ensures that the sum in the right-hand side of \eqref{eq:diffl} is well-defined. 

However, it is sometimes possible to work over smaller coefficient fields. 
One such setting
is that of {\em exact} Lagrangian submanifolds in an exact symplectic
manifold. Namely, assume that $\omega=d\theta$ for some 1-form $\theta$ on $M$, and 
there exist functions $f_i\in C^\infty(L_i,\R)$ 
such that $\theta_{|L_i}=df_i$ (for $i=0,1$). Then, by Stokes' theorem,
any strip connecting intersection points $p$
and $q$ satisfies $\int u^*\omega=(f_1(q)-f_0(q))-(f_1(p)-f_0(p))$. Thus,
rescaling each generator by $p\mapsto T^{f_1(p)-f_0(p)}p$, we can eliminate
the weights $T^{\omega([u])}$ from \eqref{eq:diffl}, and work directly over
the coefficient field $\mathbb{K}$ instead of $\Lambda$.
\end{remark}

Floer's construction \cite{Floer} was subsequently extended to more general
settings, beginning with Oh's result on monotone Lagrangians \cite{Oh}, and culminating with the sophisticated methods introduced by
Fukaya, Oh, Ohta and Ono for the general case \cite{FO3book}; however as we
will see below, Theorem \ref{thm:floer} does not hold in full generality, as
pseudo-holomorphic discs with boundary in $L_0$ or $L_1$ ``obstruct'' Floer
cohomology.

\subsection{Maslov index and grading}
The Maslov index plays a similar role in the index formula for
pseudo-holomorphic discs to that played by the first Chern class in that
for closed pseudo-holomorphic curves; in fact it can be viewed as a relative
version of the Chern class.

Denote by $LGr(n)$ the Grassmannian of Lagrangian $n$-planes in the
symplectic vector space $(\R^{2n},\omega_0)$. It is a classical fact that
the unitary group $U(n)$ acts transitively on $LGr(n)$, so that 
$LGr(n)\simeq U(n)/O(n)$, from which it follows by an easy calculation
that $\pi_1(LGr(n))\simeq \Z$ (see e.g.\ \cite{McS1}). This can be understood
concretely by using the square of the determinant map, 
$\det^2:U(n)/O(n)\to S^1$, which induces an isomorphism on fundamental
groups; the Maslov index of a loop in $LGr(n)$ is simply the winding number
of its image under this map.

In a similar vein, consider two paths $\ell_0,\ell_1:[0,1]\to LGr(n)$ of
Lagrangian subspaces in $\R^{2n}$, such that $\ell_0(0)$ is transverse to
$\ell_1(0)$ and $\ell_0(1)$ is transverse to $\ell_1(1)$.
The Maslov index of the path $\ell_1$ relative to $\ell_0$ is then the
number of times (counting with signs and multiplicities) at which 
$\ell_0(t)$ and $\ell_1(t)$ are not transverse to each other. (More precisely,
it is the intersection number of the path $(\ell_0(t),\ell_1(t))$ with the
hypersurface in $LGr(n)\times LGr(n)$ consisting of non-transverse pairs
of subspaces.)

We now return to our main discussion,
and consider a map $u:\R\times [0,1]\to M$ satisfying
the boundary conditions \eqref{eq:boundary}. Since $\R\times [0,1]$
is contractible, the pullback $u^*TM$ is a trivial symplectic vector
bundle; fixing a trivialization, we can view 
$\ell_0=u_{|\R\times\{0\}}^*TL_0$ and $\ell_1=u_{|\R\times\{1\}}^*TL_1$ 
as paths
(oriented with $s$ going from $+\infty$ to $-\infty$)
in $LGr(n)$, one connecting $T_pL_0$ to $T_qL_0$ and the
other connecting $T_pL_1$ to $T_qL_1$. The index of $u$ can then
be defined as the Maslov index of the path $\ell_1$
relative to $\ell_0$.

An equivalent definition, which generalizes more readily to the discs that appear
in the definition of product operations, is as follows. Given a pair of
transverse subspaces $\lambda_0,\lambda_1\in LGr(n)$, and identifying $\R^{2n}$
with $\C^n$, there exists an element
$A\in Sp(2n,\R)$ which maps $\lambda_0$ to $\R^n\subset \C^n$ and
$\lambda_1$ to $(i\R)^n\subset \C^n$.  The subspaces
$\lambda_t=A^{-1}((e^{-i\pi t/2}\R)^n)$, $t\in [0,1]$ then provide a
distinguished homotopy class of path connecting
$\lambda_0$ to $\lambda_1$ in $LGr(n)$, which we call the {\em canonical
short path}.

\begin{definition}\label{def:index}
Given $p,q\in L_0\cap L_1$, denote by $\lambda_p$ the canonical short path
from $T_pL_0$ to $T_pL_1$ and by $\lambda_q$ that from $T_qL_0$ to $T_qL_1$.
Given a strip $u:\R\times [0,1]\to M$ connecting $p$ to $q$, 
for $i\in\{0,1\}$, denote by $\ell_i$ the path $u^*_{|\R\times\{i\}}TL_i$
oriented with $s$ going from $+\infty$ to $-\infty$, from $T_pL_i$ to
$T_qL_i$. View all these as paths in
$LGr(n)$ by fixing a trivialization of $u^*TM$.
The {\em index} of the strip $u$ is then the Maslov index of the closed loop in $LGr(n)$ 
(based at $T_qL_0$) obtained
by concatenating the paths $-\ell_0$ (i.e.\ $\ell_0$ backwards), 
$\lambda_p$, $\ell_1$, and finally $-\lambda_q$.
\end{definition}

\begin{example}
Let $M=\R^2$, and consider the strip $u$ depicted in Figure \ref{fig:strip}:
then it is an easy exercise to check, using either definition, that $\mathrm{ind}(u)=1$.
\end{example}

We now discuss the related issue of equipping Floer complexes with a
grading. In order to obtain a $\Z$-grading on $CF(L_0,L_1)$, one needs
to make sure that the index of a strip depends only on the difference between
the degrees of the two generators it connects, rather than on its homotopy
class. This is ensured by the following two requirements:

\begin{enumerate}
\item The first Chern class of $M$ must be 2-torsion: 
$2c_1(TM)=0$. This allows one to lift the Grassmannian $LGr(TM)$ of
Lagrangian planes in $TM$ (an $LGr(n)$-bundle over $M$) to a fiberwise
universal cover $\widetilde{LGr}(TM)$, the Grassmannian of {\em graded
Lagrangian planes} in $TM$ (an $\widetilde{LGr}(n)$-bundle over $M$).

Concretely, given a nowhere vanishing section $\Theta$ of $(\Lambda^n_\C
T^*M)^{\otimes 2}$, the argument of $\Theta$ associates to any Lagrangian
plane $\ell$ a {\em phase} $\varphi(\ell)=\arg(\Theta_{|\ell})\in S^1=\R/2\pi\Z$;
a graded lift of $\ell$ is the choice of a real lift of
$\tilde\varphi(\ell)\in\R$ of $\varphi(\ell)$.\smallskip

\item The {\em Maslov class} of $L$, $\mu_L\in \mathrm{Hom}(\pi_1(L),\Z)=H^1(L,\Z)$, vanishes.
The Maslov class is by definition the obstruction to consistently choosing
graded lifts of the tangent planes to $L$, i.e.\ lifting the
section of $LGr(TM)$ over $L$ given by $p\mapsto T_pL$ to
a section of the infinite cyclic cover $\widetilde{LGr}(TM)$. 
The Lagrangian submanifold $L$ together with the choice of such a lift is
called a {\em graded Lagrangian submanifold} of $M$.

Equivalently, given a nowhere vanishing section of
$(\Lambda^n_\C T^*M)^{\otimes 2}$, we can associate to $L$ its
{\em phase function} $\varphi_L:L\to S^1$, which maps $p\in L$ to
$\varphi(T_pL)\in S^1$; the Maslov class is then the homotopy class 
$[\varphi_L]\in [L,S^1]=H^1(L,\Z)$, and a graded lift of $L$ is the
choice of a lift $\tilde\varphi_L:L\to \R$.
\end{enumerate}

When these two assumptions are satisfied, fixing graded lifts $\tilde{L}_0,\tilde{L}_1$ of the Lagrangian
submanifolds $L_0,L_1\subset M$ determines a natural $\Z$-grading on
the Floer complex $CF(L_0,L_1)$ as follows. For all $p\in L_0\cap L_1$,
we obtain a preferred homotopy class of path connecting $T_pL_0$ to
$T_pL_1$ in $LGr(T_pM)$ by connecting the 
chosen graded lifts 
of the tangent spaces at $p$ via a path in $\widetilde{LGr}(T_pM)$.
Combining this path with $-\lambda_p$ (the
canonical short path from $T_pL_0$ to $T_pL_1$, backwards),
we obtain a closed loop in $LGr(T_pM)$; the degree of $p$ is by definition
the Maslov index of this loop. It is then easy to check that any strip
connecting $p$ to $q$ satisfies
\begin{equation}
\mathrm{ind}(u)=\deg(q)-\deg(p).
\end{equation}
In particular the Floer differential \eqref{eq:diffl} has degree 1.

In general, if we do not restrict ourselves to symplectic manifolds with
torsion $c_1(TM)$ and Lagrangian submanifolds with vanishing Maslov class,
the natural grading on Floer cohomology is only by a finite cyclic group.
As an important special case, if we simply assume that the Lagrangian submanifolds
$L_0,L_1$ are oriented, then we have a $\Z/2$-grading, where the degree of
a generator $p$ of $CF(L_0,L_1)$ is determined by the sign of the intersection
between $L_0$ and $L_1$ at $p$: namely $\deg(p)=0$ if the canonical short path from
$T_pL_0$ to $T_pL_1$ maps the given orientation of $T_pL_0$ to that of
$T_pL_1$, and $\deg(p)=1$ otherwise.

(Another approach, which we won't discuss further,
is to enlarge the coefficient
field by a formal variable of non-zero degree to keep track of
the Maslov indices of different homotopy classes.
In the monotone case, where index is proportional to symplectic area, 
it suffices to give a non-zero degree to the Novikov parameter~$T$.)

\subsection{Transversality and compactness}\label{ss:transv}

We now discuss very briefly the fundamental technical issues of
transversality and compactness.

{\em Transversality} of the moduli spaces of pseudo-holomorphic strips, 
i.e.\ the surjectivity of the linearized $\dbar$ operator at all solutions,
is needed in order to ensure that the spaces 
$\widehat{\mathcal{M}}(p,q;[u],J)$ (and other moduli spaces we will 
introduce below) are smooth manifolds of the expected dimension.
Still assuming that $L_0$ and $L_1$ intersect transversely, 
transversality for strips can be achieved by replacing the 
fixed almost-complex structure $J$ in the Cauchy-Riemann equation
\eqref{eq:dbar} by a generic
{\em family} of $\omega$-compatible almost-complex structures which 
depend on the coordinate $t$ in the strip $\R\times[0,1]$.

A more basic issue is that of defining Floer cohomology for
Lagrangian submanifolds which do not intersect transversely (in particular,
one would like to be able to define the Floer cohomology of a Lagrangian with
itself, i.e.\ the case $L_0=L_1$). In view of the requirement of Hamiltonian
isotopy invariance of the construction, the simplest approach is to
introduce an inhomogeneous Hamiltonian perturbation term into the
holomorphic curve equation: we fix a generic 
Hamiltonian $H\in C^\infty([0,1]\times M,\R)$, and consider the modified
equation $(du-X_H\otimes dt)^{0,1}=0$, i.e.\ 
\begin{equation}\label{eq:floereq}
\frac{\partial u}{\partial s}+J(t,u)\left(\frac{\partial u}{\partial t}-
X_H(t,u)\right)=0,
\end{equation}
still subject to the boundary conditions $u(s,0)\in L_0$ and $u(s,1)\in
L_1$ and a finite energy condition. For $s\to \pm\infty$, the strip $u$
converges no longer to intersection points but rather to trajectories 
of the flow of $X_H$ which start on $L_0$ and end on $L_1$:
thus the generators of the Floer complex $CF(L_0,L_1)$ are in fact defined
to be flow lines $\gamma:[0,1]\to M$,
$\dot\gamma(t)=X_H(t,\gamma(t))$, such that $\gamma(0)\in L_0$ and $\gamma(1)\in L_1$.
Equivalently, by considering $\gamma(0)$, we set 
$\mathcal{X}(L_0,L_1)=
L_0\cap (\phi_H^1)^{-1}(L_1)$, where $\phi_H^1\in \mathrm{Ham}(M,\omega)$ is the time 1 flow
generated by $H$. In this sense, the generators are {\em perturbed}
intersection points of $L_0$ with $L_1$, where the perturbation is given by the
Hamiltonian diffeomorphism $\phi_H^1$.

\begin{remark}\label{rmk:hamtrick}
The perturbed equation \eqref{eq:floereq} can be recast as a plain
Cauchy-Riemann equation by the following trick: consider
$\tilde{u}(s,t)=(\phi_H^t)^{-1}(u(s,t))$, where $\phi_H^t$ is the flow of
$X_H$ over the interval $[0,t]$. Then $$\frac{\partial \tilde{u}}{\partial t}=
(\phi_H^t)^{-1}_*\left(\frac{\partial u}{\partial t}-X_H\right),$$ 
so Floer's equation \eqref{eq:floereq}
becomes
$$\frac{\partial \tilde{u}}{\partial s}+\tilde{J}(t,\tilde{u})\,\frac{\partial
\tilde{u}}{\partial t}=0,$$
where $\tilde{J}(t)=(\phi_H^t)^{-1}_*(J(t))$. Hence solutions to Floer's
equation correspond to
honest $\tilde{J}$-holomorphic strips with boundaries on $L_0$ and
$(\phi_H^1)^{-1}(L_1)$.
\end{remark}

{\em Compactness} of the moduli spaces is governed by Gromov's compactness
theorem, according to which any sequence of $J$-holomorphic curves with 
uniformly bounded energy admits a subsequence which converges, up to
reparametrization, to a nodal {\em tree} of $J$-holomorphic curves. 
The components of the limit curve are obtained as limits of
different reparametrizations of the given sequence of curves, focusing on
the different regions of the domain in which a non-zero amount of energy
concentrates (``bubbling''). 
In the case of a sequence of $J$-holomorphic
strips \hbox{$u_n:\R\times [0,1]\to M$} with boundary on Lagrangian submanifolds $L_0$ and $L_1$, there
are three types of phenomena to consider:

\begin{figure}[b]
\setlength{\unitlength}{5mm}%
\begin{picture}(14,3.8)(2,-1.1)
\newgray{gray30}{0.8}
\newgray{gray15}{0.9}
\psset{unit=\unitlength}
\pscurve[fillstyle=solid,linestyle=none,fillcolor=gray15](8.75,0.75)(10.5,2)(12.5,2)(14.25,0.75)
\pscurve[fillstyle=solid,linestyle=none,fillcolor=gray15](8.75,0.75)(10.5,-0.5)(12.5,-0.5)(14.25,0.75)
\pscurve[fillstyle=solid,linestyle=none,fillcolor=gray15](8.75,0.75)(7,2)(5,2)(3.25,0.75)
\pscurve[fillstyle=solid,linestyle=none,fillcolor=gray15](8.75,0.75)(7,-0.5)(5,-0.5)(3.25,0.75)
\pscurve(2,2)(3.25,0.75)(5,-0.5)(7,-0.5)(8.75,0.75)(10.5,2)(12.5,2)(14.25,0.75)(15.5,-0.5)
\pscurve(2,-0.5)(3.25,0.75)(5,2)(7,2)(8.75,0.75)(10.5,-0.5)(12.5,-0.5)(14.25,0.75)(15.5,2)
\pscircle*(3.25,0.75){0.1}
\pscircle*(8.75,0.75){0.1}
\pscircle*(14.25,0.75){0.1}
\put(11,-1.25){\small $L_0$}
\put(11,2.4){\small $L_1$}
\put(5.5,-1.25){\small $L_0$}
\put(5.5,2.4){\small $L_1$}
\put(3.1,1.2){\small $q$}
\put(8.7,1.1){\small $r$}
\put(14.2,1.2){\small $p$}
\end{picture} \qquad
\begin{picture}(9,4.8)(7,-1.1)
\newgray{gray30}{0.8}
\newgray{gray15}{0.9}
\psset{unit=\unitlength}
\pscurve[fillstyle=solid,linestyle=none,fillcolor=gray15](8.75,0.75)(10.5,2)(12.5,2)(14.25,0.75)
\pscurve[fillstyle=solid,linestyle=none,fillcolor=gray15](8.75,0.75)(10.5,-0.5)(12.5,-0.5)(14.25,0.75)
\pscurve(7.5,-0.5)(8.75,0.75)(10.5,2)(12.5,2)(14.25,0.75)(15.5,-0.5)
\pscurve(7.5,2)(8.75,0.75)(10.5,-0.5)(12.5,-0.5)(14.25,0.75)(15.5,2)
\pscircle*(8.75,0.75){0.1}
\pscircle*(14.25,0.75){0.1}
\pscircle[fillstyle=solid,fillcolor=gray15](12.5,3.05){1}
\put(11,-1.25){\small $L_0$}
\put(10.6,2.7){\small $L_1$}
\put(8.6,1.2){\small $q$}
\put(14.2,1.2){\small $p$}
\pscircle*(12.3,2.05){0.1}
\end{picture}
\caption{Possible limits of pseudo-holomorphic strips: a broken strip (left) and a disc bubble (right).}
\label{fig:broken}
\end{figure}
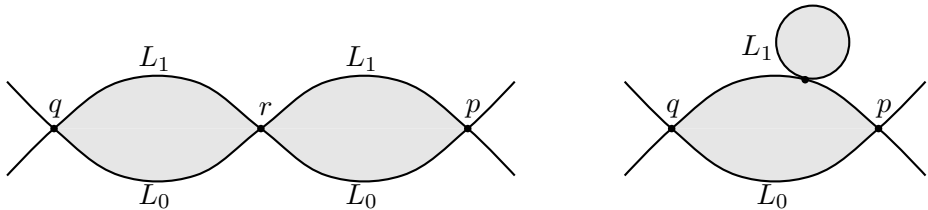

\begin{enumerate}
\item {\em strip breaking}: energy concentrates at either end $s\to \pm\infty$, i.e.\ there
is a sequence $a_n\to \pm\infty$ such that the translated strips
$u_n(s-a_n,t)$ converge to a non-constant limit strip (Figure
\ref{fig:broken} left);

\item {\em disc bubbling}: energy concentrates at a point on the boundary
of the strip ($t\in \{0,1\}$), where suitable rescalings of $u_n$ converge
to a $J$-holomorphic disc in $M$ with boundary entirely contained in either
$L_0$ or $L_1$ (Figure \ref{fig:broken} right);

\item {\em sphere bubbling}: energy concentrates at an interior point of
the strip, where suitable rescalings of $u_n$ converge to a $J$-holomorphic
sphere in $M$.
\end{enumerate}

\noindent
As we will see below, strip breaking is the key geometric ingredient in the
proof that the Floer differential squares to zero, {\em provided that
disc bubbling can be excluded}. This is not simply a technical issue -- in
general the Floer differential does not square to zero, as illustrated by 
Example \ref{ex:obstructed} below.  Another issue posed by disc and sphere 
bubbling is that of transversality: the perturbation techniques we have
outlined above are in general not sufficient to achieve transversality 
for limit curves that include disc or sphere bubble components. More 
sophisticated techniques, such as those proposed by Fukaya et al 
\cite{FO3book}\footnote{The cautious reader should be aware that, as of this writing, 
the analytic foundations of this approach are still the subject of some 
controversy.}, or the {\em polyfolds} developed by
Hofer-Wysocki-Zehnder \cite{Hofer}, are needed to extend 
Lagrangian Floer theory to the greatest possible level of generality.

In our case, the absence of disc and sphere bubbles is ensured by the 
assumption that $[\omega]\cdot \pi_2(M,L_i)=0$ in the statement 
of Theorem \ref{thm:floer}. A more general context in which the theory
still works is when bubbling can be excluded for dimension reasons,
for instance when all bubbles are guaranteed to have Maslov index greater
than $2$.
(The important limit case where the minimal Maslov index is equal to 2 can also 
be handled by elementary methods; however, in that case disc bubbling can
occur and the Floer differential does not automatically square to zero.)
A common setting where an {\em a priori} lower bound on the Maslov
index can be guaranteed is that of {\em monotone} Lagrangian submanifolds in
monotone symplectic manifolds, i.e.\ when the symplectic area of discs
and their Maslov index are proportional to each other \cite{Oh}.

\subsection{Sketch of proof of Theorem \ref{thm:floer}}

The proof that the Floer differential squares to
zero (under the assumption that disc and sphere bubbling cannot occur)
is conceptually similar to that for Morse (co)homology.

Fix Lagrangian submanifolds $L_0$ and $L_1$ as in Theorem \ref{thm:floer}, 
a generic almost-complex structure $J$ and a Hamiltonian perturbation 
$H$ so as to ensure transversality.
Given two generators $p,q$ of the Floer complex, and a homotopy class
$[u]$ with $\mathrm{ind}([u])=2$, the moduli space
$\mathcal{M}(p,q;[u],J)$ is a 1-dimensional manifold.
Since our assumptions exclude the possibilities of disc or sphere bubbling,
Gromov compactness implies that this moduli space can be compactified
to a space $\overline{\mathcal{M}}(p,q;[u],J)$ whose elements are broken
strips connecting $p$ to $q$ and representing the total class $[u]$. 

Two-component broken strips of the sort depicted
in Figure \ref{fig:broken} (left) correspond to products of moduli spaces
$\mathcal{M}(p,r;[u'],J)\times \mathcal{M}(r,q;[u''],J)$, where $r$ is any
generator of the Floer complex and $[u']+[u'']=[u]$. Observe that the index is
additive under such decompositions; moreover, transversality implies that any
non-constant strip must have index at least 1. Thus, the only possibility 
is $\mathrm{ind}([u'])=\mathrm{ind}([u''])=1$, and broken
configurations with more than two components cannot occur.

Conversely, a {\em gluing theorem} states that every broken
strip is locally the limit of a unique family of index 2 strips, and
$\overline{\mathcal{M}}(p,q;[u],J)$ is a 1-dimensional manifold with
boundary, with
\begin{equation}\label{eq:brokenM}
\partial \overline{\mathcal{M}}(p,q;[u],J) = \coprod_{\substack{
r\in \mathcal{X}(L_0, L_1)\\ [u']+[u'']=[u]\\
\!\!\mathrm{ind}([u'])=\mathrm{ind}([u''])=1\!\!}}
\Bigl(\mathcal{M}(p,r;[u'],J)\times \mathcal{M}(r,q;[u''],J)\Bigr)
\end{equation}
Moreover, the choice of orientations and spin structures on $L_0$ and $L_1$
equips all these moduli spaces with natural orientations, and
\eqref{eq:brokenM} is compatible with these orientations (up to an overall
sign). Since the total (signed) number of boundary points of a compact 
1-manifold with boundary is always zero, we conclude that
\begin{equation}\label{eq:sumcancel}
\sum_{\substack{
r\in \mathcal{X}(L_0,L_1)\\ [u']+[u'']=[u]\\
\!\!\mathrm{ind}([u'])=\mathrm{ind}([u''])=1\!\!}}
(\#\mathcal{M}(p,r;[u'],J))\,(\#\mathcal{M}(r,q;[u''],J))\,
T^{\omega([u'])+\omega([u''])}=0.
\end{equation}
Summing over all possible $[u]$, the left-hand side is precisely the coefficient of $q$ in
$\partial^2(p)$; therefore $\partial^2=0$.

When $L_0$ and/or $L_1$ bound $J$-holomorphic discs, the sum
\eqref{eq:sumcancel} no longer cancels, because the boundary of 
the 1-dimensional moduli space $\mathcal{M}(p,q;[u],J)$ also contains
configurations with disc bubbles. The following example shows that this
is an issue even in the monotone case.

\begin{example}\label{ex:obstructed}
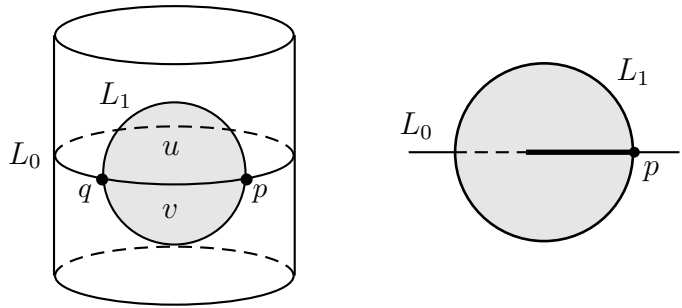
\begin{figure}[b]
\setlength{\unitlength}{8mm}
\begin{picture}(4,4.7)(-2,-2.2)
\newgray{gray30}{0.8}
\newgray{gray15}{0.9}
\psset{unit=\unitlength}
\psellipticarc(0,-2)(2,0.5){-180}{0}
\psellipticarc[linestyle=dashed](0,-2)(2,0.5){0}{180}
\pscircle[fillstyle=solid,fillcolor=gray15](0,-0.3){1.2}
\psellipticarc[linestyle=dashed](0,0)(2,0.5){0}{180}
\psellipticarc(0,0)(2,0.5){-180}{0}
\psellipse(0,2)(2,0.5)
\psline(-2,-2)(-2,2)
\psline(2,-2)(2,2)
\pscircle*(-1.2,-0.4){0.1}
\pscircle*(1.2,-0.4){0.1}
\put(-2.2,0){\makebox(0,0)[rc]{$L_0$}}
\put(-0.7,1){\makebox(0,0)[rc]{$L_1$}}
\put(-0.2,0){$u$}
\put(-0.2,-1.05){$v$}
\put(-1.6,-0.7){\small $q$}
\put(1.3,-0.7){\small $p$}
\end{picture}
\qquad\quad
\setlength{\unitlength}{1.2cm}
\begin{picture}(3,3)(-1.5,-1.5)
\psset{unit=\unitlength}
\newgray{gray30}{0.8}
\newgray{gray15}{0.9}
\psline(-1.5,0)(1.5,0)
\pscircle[fillcolor=gray15,fillstyle=solid,linewidth=1pt](0,0){1}
\psline[linestyle=dashed](-1,0)(1,0)
\psline[linewidth=2pt](-0.2,0)(1,0)
\pscircle*(1,0){0.06}
\put(1.1,-0.25){$p$}
\put(-1.6,0.2){$L_0$}
\put(0.8,0.8){$L_1$}
\end{picture}
\caption{A counterexample to $\partial^2=0$}
\label{fig:obstructed}
\end{figure}

Consider again the cylinder $M=\R\times S^1$, and let $L_0$ be a simple
closed curve that goes around the cylinder once, and $L_1$ a homotopically
trivial curve intersecting $L_0$ in two points $p$ and $q$, as shown in
Figure \ref{fig:obstructed} left. Then $L_0$ and $L_1$ bound precisely two
holomorphic strips of index 1, denoted by $u$ and $v$ in 
Figure~\ref{fig:obstructed}. (There are other holomorphic discs with boundary on
$L_0$ and $L_1$ but those have higher index.) Comparing with the convention
depicted in Figure \ref{fig:strip}, $u$ is a trajectory from $p$ to $q$,
while $v$ is a trajectory from $q$ to $p$: thus we have
$$\partial p=\pm T^{\omega(u)}\,q\quad \text{and}\quad
\partial q=\pm T^{\omega(v)}\,p,$$
and $\partial^2\neq 0$. To understand why $\partial^2(p)\neq 0$,
consider the moduli space of index 2 holomorphic strips connecting $p$ to
itself. The images of these strips exactly cover the disc bounded by $L_1$,
with a slit along $L_0$, as shown in Figure \ref{fig:obstructed} right.

We can give an explicit description using local coordinates
in which $L_0$ corresponds to the real axis and $L_1$ to the unit circle:
using the upper half-disc minus the points $\pm 1$ as domain of our
maps instead of the usual $\R\times [0,1]$ (to which it is biholomorphic),
one easily checks that any index 2 
strip connecting $p$ to itself can be parametrized as
$$u_\alpha(z)=\frac{z^2+\alpha}{1+\alpha z^2}$$
for some $\alpha\in (-1,1)$ (corresponding to the end point of the slit).

The two ends of this moduli space are different: when $\alpha\to -1$,
energy concentrates at $z=\pm 1$, and
the index 2 strips $u_\alpha$ converge to a broken strip whose nonconstant
components are the index 1 strips $u$ and $v$; whereas for $\alpha\to 1$
the maps $u_\alpha$ exhibit disc bubbling at $z=i$, the limit being a
constant strip at $p$ together with a disc bubble whose image is the
disc bounded by $L_1$. Thus, broken strips do not cancel in pairs in the
manner needed for $\partial^2=0$ to hold.
\end{example}

Once the Floer differential is shown to square to zero, it remains to
prove that Floer cohomology does not depend on the choice of almost-complex 
structure and Hamiltonian perturbation. Recall that the spaces of such
choices are contractible. Thus, given two choices $(H,J)$
and $(H',J')$ (for which we assume transversality holds), let
$(H(\tau),J(\tau))$, $\tau\in [0,1]$ be a (generically chosen) 
smooth family which agrees with $(H,J)$ for $\tau=0$ and $(H',J')$ for
$\tau=1$. One can then
construct a {\em continuation map} $F:CF(L_0,L_1;H,J)\to
CF(L_0,L_1;H',J')$ by counting solutions to the equation
\begin{equation}\label{eq:cont}
\frac{\partial u}{\partial s}+J(\tau(s),t,u)\left(\frac{\partial u}{\partial t}-
X_H(\tau(s),t,u)\right)=0,
\end{equation}
where $\tau(s)$ is a smooth function of $s$ which equals $1$ for $s\ll 0$
and $0$ for $s\gg 0$. Unlike \eqref{eq:floereq}, the equation \eqref{eq:cont} 
is not invariant under translations in the $s$ direction. Given generators
$p\in\mathcal{X}(L_0,L_1;H)$ and $p'\in\mathcal{X}(L_0,L_1;H')$ of the
respective Floer complexes, the coefficient of $p'$ in $F(p)$ is defined
as a count of index 0 solutions to \eqref{eq:cont} which converge to 
$p$ at $s\to +\infty$ and to $p'$ at $s\to -\infty$ (weighted by energy as
usual).

The proof that $F$ is a chain map, i.e.\ satisfies $\partial'\circ F=
F\circ \partial$ (again assuming the absence of bubbling), comes from
studying spaces of index 1 solutions to \eqref{eq:cont}. These spaces are
1-dimensional manifolds, whose end points correspond to broken trajectories
where the main component is an index 0 solution to \eqref{eq:cont}, either
preceded by an index 1 $J$-holomorphic strip with perturbation data $H$
(if energy concentrates at $s\to +\infty$), or followed by an index
1 $J'$-holomorphic strip with perturbation data $H'$ (if energy concentrates
at $s\to -\infty$). The composition $F\circ \partial$ counts the first type
of limit configuration, while $\partial'\circ F$ counts the second type of
limit configuration, and the equality between these two maps follows again
from the statement that the total (signed) number of end points of a compact
1-manifold with boundary is zero.

Using the reverse homotopy, i.e., considering \eqref{eq:cont} with
$\tau(s)=0$ for $s\ll 0$ and $1$ for $s\gg 0$, one similarly defines 
a chain map $F':CF(L_0,L_1;H',J')\to CF(L_0,L_1;H,J)$. The chain maps
$F$ and $F'$ are quasi-inverses, i.e.\ $F'\circ F$ is
homotopic to identity (and similarly for $F\circ F'$). An explicit 
homotopy can be obtained by counting
index $-1$ solutions to a one-parameter family of equations similar to
\eqref{eq:cont} but where $\tau(s)$ is 0 near $\pm \infty$ and is nonzero
over an interval of values of $s$ of varying width.

\subsection{The Floer cohomology $HF(L,L)$}

The Floer cohomology of a Lagrangian submanifold with itself is of particular
interest in the context of Arnold's conjecture. By Weinstein's
Lagrangian neighborhood theorem, a neighborhood of a
Lagrangian submanifold $L$ in $(M,\omega)$ is symplectomorphic to a
neighborhood of the zero section of the cotangent bundle $T^*L$ with its
standard symplectic form. In light of this, we first consider the model
case of the cotangent bundle.

\begin{example}\label{ex:cotangent}
Let $N$ be a compact real $n$-dimensional manifold, and consider the
cotangent bundle $T^*N$, with its standard exact symplectic form 
(given locally by
$\omega=\sum dq_i\wedge dp_i$, where $q_i$ are local coordinates on $N$ and
$p_i$ are the dual coordinates on the fibers of the cotangent bundle).
Let $L_0$ be the zero section, and given a Morse function $f:N\to \R$ and
a small $\epsilon>0$, denote by $L_1$ the graph of the exact 1-form
$\epsilon\,df$. Then $L_0,L_1$ are exact Lagrangian submanifolds of $T^*N$,
Hamiltonian isotopic to each other (the Hamiltonian isotopy is generated by
$H=\epsilon\,f\circ \pi$ where $\pi:T^*N\to N$ is the bundle map);
$L_0$ and $L_1$ intersect transversely at the critical points of
$f$.

Choosing a graded lift of $L_0$, and transporting it through the Hamiltonian
isotopy to define a graded lift of $L_1$, we obtain a grading on the Floer
complex $CF(L_0,L_1)$; by an explicit calculation, a critical 
point $p$ of $f$ of Morse index $i(p)$ defines a generator of the Floer 
complex of degree $\deg(p)=n-i(p)$.  Thus, the grading on the Floer complex
agrees with that on the complex $CM^*(f)$ which defines the Morse cohomology 
of $f$.

The Morse differential counts index 1 trajectories of the gradient flow 
between critical points of $f$, and depends on the choice of a Riemannian
metric $g$ on $N$, which we assume to satisfy the Morse-Smale transversality
condition. A result of Floer \cite{Floer2} is that, for a suitable choice of
(time-dependent) almost-complex structure $J$ on $T^*N$, solutions of
Floer's equation
$$\frac{\partial u}{\partial s}+J(t,u)\frac{\partial u}{\partial t}=0$$
with boundary on $L_0$ and $L_1$ are regular and in one-to-one correspondence
with gradient flow trajectories
$$\dot\gamma(s)=\epsilon \nabla f(\gamma(s))$$
on $N$, the correspondence being given by $\gamma(s)=u(s,0)$.
(Note: an ascending gradient flow line with $\gamma(s)$ converging to
$p$ as $s\to +\infty$ and $q$ as $s\to -\infty$ counts as a trajectory
from $p$ to $q$ in the Morse differential.) 

To understand this correspondence between moduli spaces,
observe that, at any point $x$ of the zero section, 
the natural almost-complex structure on $T^*N$ induced by 
the metric $g$ maps the horizontal vector $\epsilon\nabla f(x)\in
T_xN\subset T_x(T^*N)$ to the vertical vector $X_H(x)=\epsilon\,df(x)\in T^*_x
N\subset T_x(T^*N)$.  This allows us to construct particularly simple
solutions of \eqref{eq:floereq} for this almost-complex structure 
and the Hamiltonian perturbation $-H$, with both boundaries of the strip
mapping to $L_0$: for any gradient flow line $\gamma$ of $f$,
we obtain a solution of \eqref{eq:floereq} by setting $u(s,t)=\gamma(s)$.
Floer's construction of strips with boundary on $L_0$ and $L_1$ 
is equivalent to this via Remark \ref{rmk:hamtrick}.

Thus, for specific choices of perturbation data, after a rescaling of the
generators by $p\mapsto T^{\epsilon f(p)} p$, 
the Floer complex of $(L_0,L_1)$ is isomorphic to the Morse complex of~$f$, 
and the Floer cohomology $HF^*(L_0,L_1)$ is isomorphic to the Morse cohomology 
of~$f$ (with coefficients in~$\Lambda$).
Using the independence of Floer cohomology under Hamiltonian isotopies and
the isomorphism between Morse and ordinary cohomology, we conclude that
$HF^*(L_0,L_0)\simeq HF^*(L_0,L_1)\simeq H^*(L_0;\Lambda)$.

(Since we are in the exact case, by Remark \ref{rmk:novikov} one could
actually work directly over $\mathbb{K}$
rather than over Novikov coefficients.)
\end{example}

Now we consider the general case of a compact Lagrangian submanifold $L$ in a symplectic
manifold $(M,\omega)$, under the assumption that $[\omega]\cdot
\pi_2(M,L)=0$. Energy estimates then imply that, for a sufficiently small
Hamiltonian perturbation, the pseudo-holomorphic strips that determine
the Floer cohomology $HF^*(L,L)$ must all be contained in a small tubular
neighborhood of $L$, so that the calculation of Floer cohomology reduces to
Example \ref{ex:cotangent}, and we get the following result (due to
Floer in the exact case and for $\K=\Z_2$):

\begin{proposition}\label{prop:hf=h}
If $[\omega]\cdot \pi_2(M,L)=0$, then $HF^*(L,L)\simeq H^*(L;\Lambda)$.
\end{proposition}

\noindent Together with Theorem \ref{thm:floer}, this implies Arnold's conjecture
(Theorem \ref{thm:floerarnold}).

\begin{example}
Let $L$ be the zero section in $T^*S^1=\R\times S^1$ (see Figure 1 left), 
and consider the Hamiltonian perturbation depicted in the figure, which
comes from a Morse function on $L=S^1$ with a maximum at $p$ and a minimum
at $q$. Then $L$ and $\psi(L)$ bound two index 1 holomorphic strips (shaded
on the figure), both connecting $p$ to $q$, and with equal areas. 
However, the contributions of these two strips to the
Floer differential cancel out (this is obvious over $\K=\Z_2$; when
$\mathrm{char}(\K)\neq 2$ a verification of signs is needed).
Thus, $\partial p=0$, and $HF^*(L,\psi(L))\simeq
H^*(S^1)$, as expected from Proposition \ref{prop:hf=h}. 
\end{example}

Things are different when $L$ bounds pseudo-holomorphic discs, and the
Floer cohomology $HF^*(L,L)$ (when it is defined) is in general smaller 
than $H^*(L;\Lambda)$. For example, let
$L$ be a monotone Lagrangian submanifold in a 
monotone symplectic manifold, with minimal Maslov index at least 2;
this is a setting where $HF^*(L,L)$ is well defined (though no longer
$\Z$-graded), as disc bubbles 
either do not occur at all or occur in cancelling pairs.
Using again a small multiple $\epsilon f$ of
a Morse function $f$ on $L$ as Hamiltonian perturbation, the Floer 
complex differs from the Morse complex $CM^*(f)$ by the presence of
additional terms in the differential; namely there are index~1 Floer 
trajectories representing a class in $\pi_2(M,L)$ of Maslov index $k$
and connecting a critical point $p$ of Morse index $i(p)$ to a critical point 
$q$ of index $i(q)=i(p)+k-1$. This situation was studied by Oh \cite{Oh,OhSeq},
who showed that the Floer complex is filtered by index (or equivalently energy), and there 
is a {\em spectral sequence} starting with the Morse cohomology
$HM^*(f)$ (or equivalently the ordinary cohomology of $L$),
whose successive differentials account for classes of increasing Maslov
index in $\pi_2(M,L)$, and converging to the Floer cohomology $HF^*(L,L)$.

It is often easier to study honest pseudo-holomorphic discs
with boundary on $L$, rather than solutions of Floer's equation with a
Hamiltonian perturbation, or strips with boundary on $L$ and its image
under a small isotopy. This has led to the development of
alternative constructions of $HF^*(L,L)$.
For instance,
another model for the Floer cohomology of a monotone Lagrangian submanifold
is the {\em pearl complex} first introduced in \cite{Oh2} (see also
\cite{BC}). In this model, the generators of the Floer complex are
again the critical points of a Morse function $f$ on $L$, but the
differential counts ``pearly trajectories'', which arise as limits of
Floer trajectories of the sort considered above as $\epsilon\to 0$. 
Namely, a pearly trajectory between
critical points $p$ and $q$ of $f$ consists of $r\ge 0$ 
pseudo-holomorphic discs in $M$ with boundary in $L$,
connected to each other and
to $p$ and $q$ by $r+1$ gradient flow lines of $f$ in $L$. (When there are no
discs, a pearly trajectory is simply a gradient flow line between $p$ 
and~$q$.)  Yet another model, proposed by Fukaya-Oh-Ohta-Ono \cite{FO3book},
uses a chain complex where $CF(L,L)=C_*(L)$ consists of
chains in $L$, and the differential is the sum of the classical boundary map
and a map defined in terms of moduli spaces of pseudo-holomorphic discs with
boundary on $L$. This model is computationally convenient, but requires
great care in its construction to address questions such as exactly what sort of chains are considered and,
in the general (non-monotone) case, how to achieve transversality of the 
evaluation maps.

\section{Product operations}
\subsection{The product}
Let $L_0,L_1,L_2$ be three Lagrangian submanifolds of $(M,\omega)$, which we
assume intersect each other transversely and do not bound any
pseudo-holomorphic discs. We now define a product operation on their Floer
complexes, i.e.\ a map
$$CF(L_1,L_2)\otimes CF(L_0,L_1)\longrightarrow CF(L_0,L_2).$$
Given intersection points $p_1\in \mathcal{X}(L_0,L_1)$,
$p_2\in \mathcal{X}(L_1,L_2)$, and $q\in \mathcal{X}(L_0,L_2)$, the
coefficient of $q$ in $p_2\cdot p_1$ is a weighted count of
pseudo-holomorphic discs in $M$ with boundary on $L_0\cup L_1\cup L_2$ and
with corners at $p_1,p_2,q$. More precisely, let $D$ be the closed unit
disc minus three boundary points, say for instance $z_0=-1,z_1=e^{-i\pi/3},
z_2=e^{i\pi/3}$, 
and observe that a neighborhood of each puncture in $D$ is conformally equivalent
to a strip (i.e., the product of an infinite interval with $[0,1]$).

\begin{figure}[t]
\setlength{\unitlength}{8mm}%
\begin{picture}(16.5,4)(0,-2)
\newgray{gray30}{0.8}
\newgray{gray15}{0.9}
\psset{unit=\unitlength}
\pscurve[fillstyle=solid,fillcolor=gray15,linestyle=none]%
(0,0.5)(1,0.55)(3,1.45)(4,1.5)(4,1.5)(4,0.55)%
(4,0.55)(3,0.45)(2.5,0)(3,-0.45)(4,-0.55)%
(4,-0.55)(4,-1.5)(4,-1.5)(3,-1.45)(1,-0.55)(0,-0.5)
\pscurve(0,0.5)(1,0.55)(3,1.45)(4,1.5)
\pscurve(0,-0.5)(1,-0.55)(3,-1.45)(4,-1.5)
\pscurve(4,0.55)(3,0.45)(2.5,0)(3,-0.45)(4,-0.55)
\put(1.5,-1.4){$L_0$}
\put(1.5,1.3){$L_2$}
\put(3,-0.2){$L_1$}
\pscircle[fillstyle=solid,fillcolor=gray15](7.5,0){1.5}
\pscircle*(6,0){0.1}
\pscircle*(8.25,1.3){0.1}
\pscircle*(8.25,-1.3){0.1}
\put(6.2,-1.6){\small $L_0$}
\put(6.2,1.45){\small $L_2$}
\put(9.1,-0.3){\small $L_1$}
\put(5.5,-0.3){$z_0$}
\put(8.3,-1.6){$z_1$}
\put(8.3,1.5){$z_2$}
\put(4.8,0){$\simeq$}
\psline{->}(10,0)(11.5,0)
\put(10.55,0.15){$u$}
\pspolygon[fillstyle=solid,fillcolor=gray15,linestyle=none](13,0)(15.6,-1.5)(15.6,1.5)
\psarc[fillstyle=solid,fillcolor=gray15](13,0){3}{-30}{30}
\psarc[fillstyle=solid,fillcolor=gray15](15.6,-1.5){3}{90}{150}
\psarc[fillstyle=solid,fillcolor=gray15](15.6,1.5){3}{-150}{-90}
\psarc(13,0){3}{-45}{45}
\psarc(15.6,-1.5){3}{75}{165}
\psarc(15.6,1.5){3}{-165}{-75}
\pscircle*(13,0){0.1}
\pscircle*(15.6,1.5){0.1}
\pscircle*(15.6,-1.5){0.1}
\put(13.6,1.4){$L_2$}
\put(13.6,-1.55){$L_0$}
\put(16.2,-0.2){$L_1$}
\put(12.5,-0.1){$q$}
\put(15.7,-1.85){$p_1$}
\put(15.7,1.85){$p_2$}
\end{picture}
\caption{A pseudo-holomorphic disc contributing to the product map.}
\label{fig:product}
\end{figure}
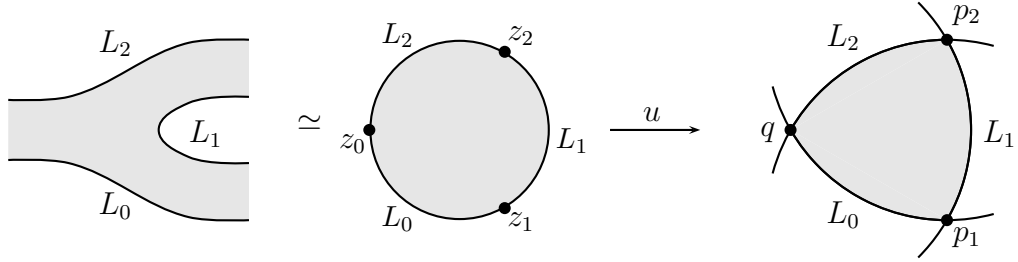

Given an almost-complex structure $J$ on $M$ and a homotopy class $[u]$, we
denote by $\mathcal{M}(p_1,p_2,q;[u],J)$ the space
of finite energy $J$-holomorphic maps $u:D\to M$ which extend continuously
to the closed disc, mapping the boundary
arcs from $z_0$ to $z_1$, $z_1$ to $z_2$, $z_2$ to $z_0$ to
$L_0,L_1,L_2$ respectively, and
the boundary punctures $z_1,z_2,z_0$ to $p_1,p_2,q$ respectively,
in the given homotopy class $[u]$ (see Figure \ref{fig:product}).

As in the case of strips, the expected dimension of
$\mathcal{M}(p_1,p_2,q;[u],J)$ is given by the index of the linearized
Cauchy-Riemann operator $D_{\dbar_J,u}$. This index can be expressed
in terms of the Maslov index, exactly as in Definition \ref{def:index}:
we now concatenate the paths given by the tangent spaces to $L_0,L_1,L_2$
going counterclockwise along the boundary of $u$, together with the appropriate canonical short
paths at $p_1,p_2,q$, to obtain a closed loop in $LGr(n)$ whose Maslov index
is equal to $\mathrm{ind}(u)$.
If $c_1(TM)$ is 2-torsion and the Maslov classes of $L_0,L_1,L_2$ vanish,
then after choosing graded lifts of the Lagrangians we have $\Z$-gradings 
on the Floer complexes, and one checks that
\begin{equation}\label{eq:indproduct}
\mathrm{ind}(u)=\deg(q)-\deg(p_1)-\deg(p_2).
\end{equation}

\begin{remark}
The apparent lack of symmetry in the index formula \eqref{eq:indproduct} 
is due to the difference between the gradings on $CF(L_0,L_2)$ and $CF(L_2,L_0)$. Namely, the given
intersection point $q\in L_0\cap L_2$ defines generators of both 
complexes, whose degrees sum to $n$ (the dimension of $L_i$). In fact, the Floer
complexes $CF(L_0,L_2)$ and $CF(L_2,L_0)$ and the differentials on them
are dual to each other, provided that the almost-complex structures and
perturbations are chosen suitably. For instance, the strip
depicted in Figure \ref{fig:strip} is a trajectory from $p$ to $q$ in the
Floer complex $CF(L_0,L_1)$,
and from $q$ to $p$ in $CF(L_1,L_0)$.
\end{remark}

Assume that transversality holds, so that the moduli spaces 
$\mathcal{M}(p_1,p_2,q;[u],J)$ are smooth manifolds; if
$\mathrm{char}(\K)\neq 2$, assume moreover that orientations and spin
structures on $L_0,L_1,L_2$ have been chosen, so as to determine orientations of the moduli
spaces. Then we define:

\begin{definition}\label{def:product}
The Floer product is the $\Lambda$-linear map $CF(L_1,L_2)\otimes
CF(L_0,L_1)\to CF(L_0,L_2)$ defined by
\begin{equation}\label{eq:product}
p_2\cdot p_1=\sum_{\substack{q\in \mathcal{X}(L_0,L_2)\\
[u]:\mathrm{ind}([u])=0}}
(\#\mathcal{M}(p_1,p_2,q;[u],J)) T^{\omega([u])}\,q.
\end{equation}
\end{definition}

\noindent As in the previous section, in general this construction needs to be modified by introducing
domain-dependent almost-complex structures and Hamiltonian perturbations
to achieve transversality. We discuss this below, but for now we assume
transversality holds without further perturbations and examine
the properties of the Floer product.

\begin{proposition}\label{prop:product}
If $[\omega]\cdot \pi_2(M,L_i)=0$ for all $i$, then the Floer product satisfies
the Leibniz rule (with suitable signs) with respect to the Floer
differentials, 
\begin{equation}\label{eq:leibniz}
\partial(p_2\cdot p_1)=\pm (\partial p_2)\cdot p_1\pm p_2\cdot (\partial
p_1),\end{equation}
and hence induces a well-defined product
$HF(L_1,L_2)\otimes HF(L_0,L_1)\to HF(L_0,L_2)$.
Moreover, this induced product on Floer cohomology groups is independent
of the chosen almost-complex structure (and Hamiltonian perturbations) and 
associative.
\end{proposition}

\noindent (However, the chain-level product on Floer complexes is {\em not}
associative, as we will see below.)

We now sketch the geometric argument behind the Leibniz rule, which relies
on an examination of index 1 moduli spaces of $J$-holomorphic discs and
their compactification. Namely, consider a triple of generators $p_1,p_2,q$
as above, and let $[u]$ be a homotopy class with $\mathrm{ind}([u])=1$.
Then (still assuming transversality) $\mathcal{M}(p_1,p_2,q;[u],J)$ is a 
smooth 1-dimensional manifold, and by Gromov compactness admits a
compactification $\overline{\mathcal{M}}(p_1,p_2,q;[u],J)$ obtained by
adding nodal trees of $J$-holomorphic curves. 

\begin{figure}[t]
\setlength{\unitlength}{7mm}%
\begin{picture}(6,5)(10.5,-2.5)
\newgray{gray30}{0.8}
\newgray{gray15}{0.9}
\psset{unit=\unitlength}
\pspolygon[fillstyle=solid,fillcolor=gray15,linestyle=none](13,0)(15.6,-1.5)(15.6,1.5)
\psarc[fillstyle=solid,fillcolor=gray15](13,0){3}{-30}{30}
\psarc[fillstyle=solid,fillcolor=gray15](15.6,-1.5){3}{90}{150}
\psarc[fillstyle=solid,fillcolor=gray15](15.6,1.5){3}{-150}{-90}
\psarc(13,0){3}{-45}{45}
\psarc(15.6,-1.5){3}{75}{150}
\psarc(15.6,1.5){3}{-150}{-75}
\psarc[fillstyle=solid,fillcolor=gray15](12,1){1.414}{-135}{-45}
\psarc[fillstyle=solid,fillcolor=gray15](12,-1){1.414}{45}{135}
\psarc(12,1){1.414}{-150}{-45}
\psarc(12,-1){1.414}{45}{150}
\pscircle*(11,0){0.1}
\pscircle*(13,0){0.1}
\pscircle*(15.6,1.5){0.1}
\pscircle*(15.6,-1.5){0.1}
\put(13.6,1.4){\small $L_2$}
\put(13.6,-1.55){\small $L_0$}
\put(11.6,0.65){\small $L_2$}
\put(11.6,-0.85){\small $L_0$}
\put(16.2,-0.2){\small $L_1$}
\put(10.5,-0.1){\small $q$}
\put(15.7,-1.85){\small $p_1$}
\put(15.7,1.85){\small $p_2$}
\end{picture}
\qquad
\begin{picture}(5,5)(12.5,-2)
\newgray{gray30}{0.8}
\newgray{gray15}{0.9}
\psset{unit=\unitlength}
\pspolygon[fillstyle=solid,fillcolor=gray15,linestyle=none](13,0)(15.6,-1.5)(15.6,1.5)
\psarc[fillstyle=solid,fillcolor=gray15](13,0){3}{-30}{30}
\psarc[fillstyle=solid,fillcolor=gray15](15.6,-1.5){3}{90}{150}
\psarc[fillstyle=solid,fillcolor=gray15](15.6,1.5){3}{-150}{-90}
\psarc(13,0){3}{-45}{30}
\psarc(15.6,-1.5){3}{90}{165}
\psarc(15.6,1.5){3}{-165}{-75}
\psarc[fillstyle=solid,fillcolor=gray15](15.6,2.9){1.4}{-90}{0}
\psarc[fillstyle=solid,fillcolor=gray15](17,1.5){1.4}{90}{180}
\pscircle*(13,0){0.1}
\pscircle*(15.6,1.5){0.1}
\pscircle*(17,2.9){0.1}
\pscircle*(15.6,-1.5){0.1}
\put(13.6,1.4){\small $L_2$}
\put(15.1,2.3){\small $L_2$}
\put(13.6,-1.55){\small $L_0$}
\put(16.2,-0.2){\small $L_1$}
\put(16.8,1.8){\small $L_1$}
\put(12.5,-0.1){\small $q$}
\put(15.7,-1.85){\small $p_1$}
\put(17.2,2.9){\small $p_2$}
\end{picture}
\qquad
\begin{picture}(5,5)(12.5,-3)
\newgray{gray30}{0.8}
\newgray{gray15}{0.9}
\psset{unit=\unitlength}
\pspolygon[fillstyle=solid,fillcolor=gray15,linestyle=none](13,0)(15.6,-1.5)(15.6,1.5)
\psarc[fillstyle=solid,fillcolor=gray15](13,0){3}{-30}{30}
\psarc[fillstyle=solid,fillcolor=gray15](15.6,-1.5){3}{90}{150}
\psarc[fillstyle=solid,fillcolor=gray15](15.6,1.5){3}{-150}{-90}
\psarc(13,0){3}{-30}{45}
\psarc(15.6,-1.5){3}{75}{165}
\psarc(15.6,1.5){3}{-165}{-90}
\psarc[fillstyle=solid,fillcolor=gray15](15.6,-2.9){1.4}{0}{90}
\psarc[fillstyle=solid,fillcolor=gray15](17,-1.5){1.4}{180}{270}
\pscircle*(13,0){0.1}
\pscircle*(15.6,1.5){0.1}
\pscircle*(17,-2.9){0.1}
\pscircle*(15.6,-1.5){0.1}
\put(13.6,1.4){\small $L_2$}
\put(15.1,-2.3){\small $L_0$}
\put(13.6,-1.55){\small $L_0$}
\put(16.2,-0.2){\small $L_1$}
\put(16.8,-1.8){\small $L_1$}
\put(12.5,-0.1){\small $q$}
\put(15.7,1.85){\small $p_2$}
\put(17.2,-2.9){\small $p_1$}
\end{picture}
\caption{The ends of a 1-dimensional moduli space $\mathcal{M}(p_1,p_2,q;[u],J)$.}
\label{fig:breaking}
\end{figure}
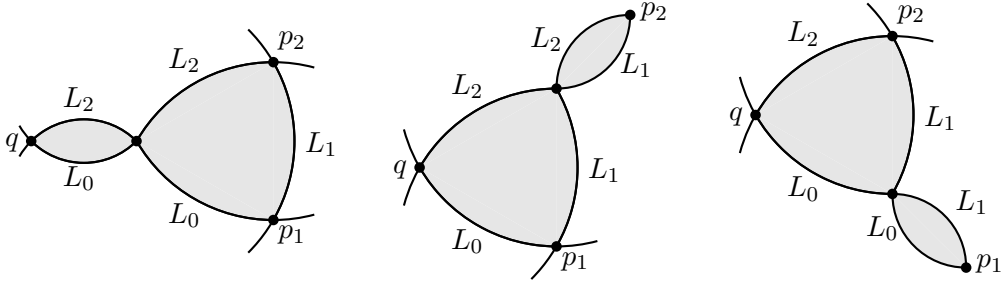

Since our assumptions exclude
bubbling of discs or spheres, the only phenomenon that can occur is strip-breaking
(when energy concentrates at one of the three ends of the punctured disc $D$).
Since transversality excludes the presence of discs of index less than 0
and nonconstant strips of index less than 1, and since the sum of the
indices of the limit components must be 1, there are only
three types of limit configurations to be considered, all consisting of
an index 0 disc with boundary on $L_0,L_1,L_2$ and an index 1 strip with
boundary on two of these three submanifolds; see Figure
\ref{fig:breaking}.

The three types of configurations contribute to the
coefficient of $T^{\omega([u])}q$ in \hbox{$\partial(p_2\cdot p_1)$}
(Figure \ref{fig:breaking} left), $(\partial p_2)\cdot p_1$
(middle), and $p_2\cdot (\partial p_1)$ (right) respectively.
On the other hand, a gluing theorem states that every such configuration
arises as an end of $\mathcal{M}(p_1,p_2,q;[u],J)$, and that the
compactified moduli space is a 1-dimensional compact manifold with boundary.
Moreover, the orientations agree up to overall sign factors depending only
on the degrees of $p_1$ and $p_2$.
Since the (signed) total number of boundary points of
$\overline{\mathcal{M}}(p_1,p_2,q;[u],J)$ is zero,
the Leibniz rule \eqref{eq:leibniz} follows.

Before moving on to higher products, we briefly discuss the issue of
transversality and compatibility in the choice of perturbations.
As in the case of strips, even without assuming that $L_0,L_1,L_2$ intersect
transversely, we can ensure transversality by introducing domain-dependent
almost-complex structures and Hamiltonian perturbations; however, for the
Leibniz rule to hold, these need to be chosen suitably
near the punctures $z_0,z_1,z_2$.
Fix once and for all ``strip-like ends'' near the punctures, i.e.\
biholomorphisms from $\R_+\times [0,1]$ (resp.\ $\R_-\times [0,1]$) 
to neighborhoods of the punctures $z_1$ and $z_2$ (resp.\ $z_0$) in $D$;
we denote by $s+it$ the natural complex coordinate in each strip-like end.
Also fix a 1-form $\beta\in \Omega^1(D)$, such that $\beta_{|\partial D}=0$ and
$\beta=dt$ in each strip-like end.
Now, given $L_0,L_1,L_2$, we choose a family of $\omega$-compatible 
almost-complex structures depending smoothly on $z\in D$, i.e.\ $J\in
C^\infty(D,\mathcal{J}(M,\omega))$, and a family of
Hamiltonians $H\in C^\infty(D\times M,\R)$, with the property that
in each strip-like end $J(z)$ and $H(z)$ depend only on the coordinate 
$t\in [0,1]$. We then perturb the Cauchy-Riemann equation to
\begin{equation}\label{eq:perthol}
\Bigl(du-X_H\otimes \beta\Bigr)^{0,1}_J=0,
\end{equation}
which in each strip-like end reduces to \eqref{eq:floereq}.

For $0\le i<j\le 2$, denote by $H_{ij}\in C^\infty([0,1]\times M,\R)$
and $J_{ij}\in C^\infty([0,1],\mathcal{J}(M,\omega))$ the time-dependent
Hamiltonians and almost-complex structures on the strip-like end whose
boundaries map to $L_i$ and $L_j$. 
The solutions of \eqref{eq:perthol} converge no longer to intersection
points of $L_i\cap L_j$, but to trajectories of the time 1 flow generated
by $H_{ij}$ which begin on $L_i$ and end on $L_j$, i.e.\ generators of
the perturbed Floer complex of $(L_i,L_j)$ with respect to the Hamiltonian
perturbation $H_{ij}$. Moreover, when strip breaking occurs, the main
component remains a solution of \eqref{eq:perthol}, while the strip
component that breaks off is a solution of \eqref{eq:floereq} with respect
to $H_{ij}$ and $J_{ij}$. 

Thus, by considering the moduli spaces of solutions to the perturbed equation 
\eqref{eq:perthol} and proceeding as in Definition \ref{def:product},
we obtain a product map 
$$CF(L_1,L_2;H_{12},J_{12})\otimes CF(L_0,L_1;H_{01},J_{01})\longrightarrow
CF(L_0,L_2;H_{02},J_{02})$$
on the perturbed Floer complexes, and Proposition \ref{prop:product} still
holds (with respect to the perturbed Floer differentials).

\subsection{Higher operations}
Given $k+1$ Lagrangian submanifolds $L_0,\dots,L_k$, a construction similar
to those above allows us to define an
operation
$$\mu^k:CF(L_{k-1},L_k)\otimes \dots\otimes CF(L_1,L_2)\otimes CF(L_0,L_1)
\longrightarrow CF(L_0,L_k)$$
(of degree $2-k$ in the situation where the Floer complexes are graded),
where $\mu^1$ is the Floer differential and $\mu^2$ is the product.

Given generators $p_i\in \mathcal{X}(L_{i-1},L_i)$ ($i=1,\dots,k$) and
$q\in \mathcal{X}(L_0,L_k)$, the coefficient of $q$ in
$\mu^k(p_k,\dots,p_1)$ is a count (weighted by area) of (perturbed)
pseudo-holomorphic discs in $M$ with boundary on $L_0\cup \dots\cup L_k$
and corners at $p_1,\dots,p_k,q$.

Specifically, one considers maps $u:D\to M$ whose domain $D$ is the 
closed unit disc minus $k+1$ boundary points $z_0,z_1,\dots,z_k\in S^1$,
lying in that order along the unit circle. The positions of these marked
points are not fixed, and the moduli space $\mathcal{M}_{0,k+1}$ of 
conformal structures on the domain $D$, i.e., the quotient of the
space of ordered $(k+1)$-tuples of points on $S^1$ by the action of
$\mathrm{Aut}(D^2)$, is a contractible $(k-2)$-dimensional manifold.

Given an almost-complex structure $J$ on $M$ and a homotopy class $[u]$,
we denote by $\mathcal{M}(p_1,\dots,p_k,q;[u],J)$ the space of
$J$-holomorphic maps $u:D\to M$ (where the positions of $z_0,\dots,z_k$
are not fixed a priori) which extend continuously to the closed disc,
mapping the boundary arcs from $z_i$ to $z_{i+1}$ (or $z_0$ for $i=k$) to 
$L_i$, and the boundary punctures $z_1,\dots,z_k,z_0$ to $p_1,\dots,p_k,q$
respectively, in the given homotopy class $[u]$, up to the action of
$\mathrm{Aut}(D^2)$ by reparametrization. (Or, equivalently, one can avoid
quotienting and instead take a slice for the reparametrization action by 
fixing the positions of three of the $z_i$.)

For a fixed conformal structure on $D$, the index of the linearized
Cauchy-Riemann operator is again given by the Maslov index, as previously.
Thus, accounting for deformations of the conformal structure on $D$,
assuming transversality, the expected dimension of the moduli space is
\begin{equation}\label{eq:dimmk}
\dim\mathcal{M}(p_1,\dots,p_k,q;[u],J)=k-2+\mathrm{ind}([u])=
k-2+\deg(q)-\textstyle\sum\limits_{i=1}^k \deg(p_i).
\end{equation}

Thus, assuming transversality, and choosing orientations and spin structures
on $L_0,\dots,L_k$ if $\mathrm{char}(\K)\neq 2$, we define:

\begin{definition}\label{def:mk}
The operation $\mu^k:CF(L_{k-1},L_k)\otimes \dots\otimes CF(L_0,L_1)
\to CF(L_0,L_k)$ is 
the $\Lambda$-linear map defined by
\begin{equation}\label{eq:mk}
\mu^k(p_k,\dots,p_1)=\sum_{\substack{q\in \mathcal{X}(L_0,L_k)\\
[u]:\mathrm{ind}([u])=2-k}}
(\#\mathcal{M}(p_1,\dots,p_k,q;[u],J)) T^{\omega([u])}\,q.
\end{equation}
\end{definition}

\begin{remark}\label{rmk:perturb}
As before, in general this construction needs
to be modified by introducing domain-dependent almost-complex structures and 
Hamiltonian perturbations to achieve transversality. Thus, we actually
count solutions of a perturbed Cauchy-Riemann equation similar to
\eqref{eq:perthol}, involving a domain-dependent almost-complex structure
$J\in C^\infty(D,\mathcal{J}(M,\omega))$ and Hamiltonian $H\in
C^\infty(D\times M,\R)$. As before, compatibility with strip-breaking
requires that, in each of the $k+1$ strip-like ends near the
punctures of $D$, the chosen $J$ and $H$ depend only on the coordinate $t\in [0,1]$ and agree with
the almost-complex structures and Hamiltonians used to construct the Floer
complexes $CF(L_i,L_{i+1})$ and $CF(L_0,L_k)$. An additional compatibility
condition comes from the possible degenerations of the domain $D$ to unions
of discs with fewer punctures, as discussed below: we need to require
that, when $D$ degenerates in such a way, $H$ and $J$ are
translation-invariant in the strip-like regions connecting
the components and agree with the choices made in the construction of the
Floer complexes $CF(L_i,L_j)$, while in each component $H$ and $J$ agree
with the choices made for that moduli space of discs with fewer punctures. 
This forces the choices of $H$ and $J$ to
further depend on the conformal structure of $D$. We refer the reader to
\cite{SeBook} for a detailed construction (and proof of existence) of
compatible and consistent choices of perturbation data $(H,J)$.
\end{remark}

The algebraic properties of $\mu^k$ follow from the study of the limit
configurations that arise in compactifications of 1-dimensional moduli
spaces of (perturbed) pseudo-holomorphic discs; besides strip breaking, 
there are now other possibilities, corresponding to cases where the domain
$D$ degenerates. The moduli space of conformal structures
$\mathcal{M}_{0,k+1}$ admits a natural compactification to a
$(k-2)$-dimensional polytope
$\overline{\mathcal{M}}_{0,k+1}$, the Stasheff {\em associahedron},
whose top-dimensional facets correspond to nodal degenerations of $D$ to
a pair of discs $D_1\cup D_2$, with each component carrying at least two
of the marked points $z_0,\dots,z_k$; and the
higher codimension faces correspond to nodal degenerations with more
components. 
\begin{figure}[b]
\setlength{\unitlength}{8mm}
\newgray{gray30}{0.8}
\newgray{gray15}{0.9}
\begin{picture}(12,4.8)(-1,-2.8)
\psset{unit=\unitlength}
\psline{|-|}(0,-2.8)(10,-2.8)
\pscircle[fillstyle=solid,fillcolor=gray15](-1,0){1}
\pscircle[fillstyle=solid,fillcolor=gray15](1,0){1}
\pscircle[fillstyle=solid,fillcolor=gray15](5.5,0){1.5}
\pscircle[fillstyle=solid,fillcolor=gray15](10,1){1}
\pscircle[fillstyle=solid,fillcolor=gray15](10,-1){1}
\psline{->}(0,-1.5)(0,-2.5)
\psline{->}(5.5,-2)(5.5,-2.5)
\psline{->}(10,-2.2)(10,-2.5)
\pscircle*(-1.5,0.866){0.1}
\pscircle*(1.5,0.866){0.1}
\pscircle*(-1.5,-0.866){0.1}
\pscircle*(1.5,-0.866){0.1}
\put(-1.85,-1.1){\makebox(0,0)[cc]{$z_0$}}
\put(1.85,-1.1){\makebox(0,0)[cc]{$z_1$}}
\put(1.85,1.1){\makebox(0,0)[cc]{$z_2$}}
\put(-1.85,1.1){\makebox(0,0)[cc]{$z_3$}}
\pscircle*(4.3,-0.9){0.1}
\pscircle*(4.3,0.9){0.1}
\pscircle*(6.7,-0.9){0.1}
\pscircle*(6.7,0.9){0.1}
\put(4,-1.25){\makebox(0,0)[cc]{$z_0$}}
\put(7,-1.25){\makebox(0,0)[cc]{$z_1$}}
\put(7,1.25){\makebox(0,0)[cc]{$z_2$}}
\put(4,1.25){\makebox(0,0)[cc]{$z_3$}}
\pscircle*(10.866,1.5){0.1}
\pscircle*(10.866,-1.5){0.1}
\pscircle*(9.134,-1.5){0.1}
\pscircle*(9.134,1.5){0.1}
\put(8.8,-1.7){\makebox(0,0)[cc]{$z_0$}}
\put(11.2,-1.7){\makebox(0,0)[cc]{$z_1$}}
\put(11.2,1.75){\makebox(0,0)[cc]{$z_2$}}
\put(8.8,1.75){\makebox(0,0)[cc]{$z_3$}}
\end{picture}
\caption{The 1-dimensional associahedron $\overline{\mathcal{M}}_{0,4}$.}
\label{fig:associahedron}
\end{figure}
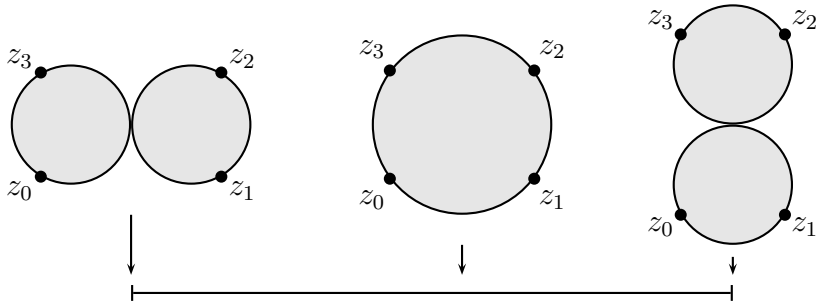

\begin{example}
$\overline{\mathcal{M}}_{0,4}$
is homeomorphic to a closed interval, whose end points correspond to
configurations where two adjacent marked points come together (Figure
\ref{fig:associahedron}). For example, fixing the positions of $z_0,z_1,z_2$
on the unit circle and letting $z_3$ vary along the arc from $z_2$ to $z_0$,
the right end point corresponds to the case where $z_3$ 
approaches $z_2$; the ``main'' component of the limit configuration carries
the marked points $z_0$ and $z_1$, while the component carrying $z_2$ and
$z_3$ arises from rescaling by suitable automorphisms of the disc.
Equivalently up to automorphisms of the disc, one could instead fix the positions of $z_1,z_2,z_3$, and let
$z_0$ vary along the arc from $z_3$ to $z_1$; the right end point then
corresponds to the case where $z_0$ approaches $z_1$.
\end{example}

\begin{proposition}\label{prop:ainfty}
If $[\omega]\cdot \pi_2(M,L_i)=0$ for all $i$, then the operations $\mu^k$
satisfy the {\em $A_\infty$-relations}
\begin{equation}\label{eq:ainfty}
\sum_{\ell=1}^{k}\sum_{j=0}^{k-\ell} (-1)^* \mu^{k+1-\ell}(p_k,\dots,
p_{j+\ell+1},\mu^\ell(p_{j+\ell},\dots,p_{j+1}),p_j,\dots,p_1)=0,
\end{equation}
where $*=j+\deg(p_1)+\dots+\deg(p_j).$
\end{proposition}

\noindent The case $k=1$ of \eqref{eq:ainfty} is the identity $\partial^2=0$, 
while $k=2$ corresponds to the Leibniz rule \eqref{eq:leibniz}. For 
$k=3$, it expresses the fact that the Floer product $\mu^2$ is associative
up to an explicit homotopy given by $\mu^3$:
\begin{multline}\label{eq:mu3} \pm (p_3\cdot p_2)\cdot p_1\pm p_3\cdot (p_2\cdot p_1)=\\
\pm \partial\mu^3(p_3,p_2,p_1)\pm \mu^3(\partial p_3,p_2,p_1)\pm
\mu^3(p_3,\partial p_2,p_1)\pm \mu^3(p_3,p_2,\partial p_1).
\end{multline}
More generally, each operation $\mu^k$ gives an explicit homotopy 
for a certain compatibility property among the preceding ones.

The proof of Proposition \ref{prop:ainfty} again relies on an analysis
of $1$-dimensional moduli spaces of (perturbed) $J$-holomorphic discs and
their compactification. Fix generators $p_1,\dots,p_k,q$ and a homotopy
class $[u]$ with $\mathrm{ind}([u])=3-k$, and assume that $J$ and $H$ 
are chosen generically (so as to achieve transversality) and compatibly
(see Remark \ref{rmk:perturb}). 
Then the moduli space $\mathcal{M}(p_1,\dots,p_k,q;[u],J)$ compactifies to
a 1-dimensional manifold with boundary, whose
boundary points correspond either to an index 1 (perturbed) $J$-holomorphic strip 
breaking off at one of the $k+1$ marked points, or to a degeneration of 
the domain to the boundary of $\overline{\mathcal{M}}_{0,k+1}$, i.e.\ to a 
pair of discs with each component carrying at least two of the marked 
points. The first case corresponds to the terms 
involving $\mu^1$ in \eqref{eq:ainfty},
while the second case corresponds to the other terms.

\begin{example}
For $k=3$, limit configurations consisting of an index~1 
strip together with an index $-1$ disc with 4 marked points account for
the right-hand side in \eqref{eq:mu3}, while those consisting of a pair
of index 0 discs with 3 marked points (when the domain degenerates to one of the two end points
of $\overline{\mathcal{M}}_{0,4}$, see Figure \ref{fig:associahedron}) account for
the two terms in the left-hand side.
\end{example}

\subsection{The Fukaya category}
There are several variants of the Fukaya category of a symplectic manifold,
depending on the desired level of generality and a number of implementation 
details. The common features are the following. The objects of the Fukaya
category are suitable Lagrangian submanifolds, equipped with extra data, 
and morphism spaces are given by Floer complexes, endowed with the Floer
differential. Composition of morphisms is given by the Floer product,
which is only associative up to homotopy, and
the Fukaya category is an {\em $A_\infty$-category}, i.e.\ the
differential and composition are the first two in a sequence of operations
$$\mu^k:\hom(L_{k-1},L_k)\otimes\dots\otimes \hom(L_0,L_1)\to
\hom(L_0,L_k)$$ (of degree $2-k$ when a $\Z$-grading is available),
satisfying the $A_\infty$-relations \eqref{eq:ainfty}.

Given the setting in which we have developed Floer theory in the preceding
sections, the most natural definition is the following:

\begin{definition}\label{def:fukaya}
Let $(M,\omega)$ be a symplectic manifold with $2c_1(TM)=0$. The objects
of the (compact) Fukaya category $\mathcal{F}(M,\omega)$ are 
compact closed, oriented, spin Lagrangian submanifolds $L\subset M$ 
such that $[\omega]\cdot \pi_2(M,L)=0$ and with vanishing Maslov class
$\mu_L=0\in H^1(L,\Z)$, together with extra data, namely the choice of
a spin structure and a graded lift of $L$. (We will usually omit those 
from the notation and simply denote the object by $L$.)

For every pair of objects $(L,L')$ (not necessarily distinct), we 
choose perturbation
data $H_{L,L'}\in C^\infty([0,1]\times M,\R)$ and $J_{L,L'}\in
C^\infty([0,1],\mathcal{J}(M,\omega))$; and for all tuples of objects
$(L_0,\dots,L_k)$ and all moduli spaces of discs, we choose consistent 
perturbation data $(H,J)$ compatible with the choices made for the
pairs of objects $(L_i,L_j)$, so as to achieve transversality for all
moduli spaces of perturbed $J$-holomorphic discs. (See \cite[\S 9]{SeBook} for
the existence of such perturbation data.)

Given this, we set $\hom(L,L')=CF(L,L';H_{L,L'},J_{L,L'})$; and the differential $\mu^1$, 
composition $\mu^2$, and higher operations $\mu^k$ are given by counts of
perturbed pseudo-holomorphic discs as in Definition \ref{def:mk}.
By Proposition \ref{prop:ainfty}, this makes $\F(M,\omega)$ a
$\Lambda$-linear, $\Z$-graded, non-unital \/{\rm (}but cohomologically unital
\cite{SeBook}{\rm )} $A_\infty$-category.
\end{definition}

\noindent One can also consider other settings: for example, we can drop the
requirement that $2c_1(TM)=0$ and the assumption of vanishing of the Maslov
class if we are content with a $\Z/2$-grading; spin structures can be
ignored if we work over a field of characteristic 2; and Novikov
coefficients are unnecessary if we restrict ourselves to exact Lagrangian
submanifolds in an exact symplectic manifold.

As is obvious from the definition, the actual chain-level details of the
Fukaya category depend very much on the choice of perturbation data;
however, the $A_\infty$-categories obtained from various choices of 
perturbation data are {\em quasi-equivalent} (i.e., they are related
by $A_\infty$-functors which induce equivalences, in fact in this case isomorphisms,
at the level of cohomology) \cite{SeBook}.

We finish this section with a few remarks.

\begin{remark}
One can recover an honest category from an $A_\infty$-category by
taking the cohomology of morphism spaces with respect to the differential
$\mu^1$; the $A_\infty$-relations imply that $\mu^2$ descends to an
associative composition operation on cohomology. 
The cohomology category of $\mathcal{F}(M,\omega)$, where 
$\hom(L,L')=HF(L,L')$ and composition is given by the 
cohomology-level Floer product, is sometimes called the Donaldson-Fukaya
category.  However, the higher operations contain important information 
that gets lost when passing to the cohomology category, and it is usually
much better to work with the chain-level $A_\infty$-category (see for
instance the next section).
\end{remark}

\begin{remark}
In the context of homological mirror symmetry, one is naturally led to
consider a slightly richer version of the Fukaya category, whose objects
are Lagrangian submanifolds equipped with {\em local systems}, i.e.\
flat vector bundles $\mathcal{E}\to L$ with unitary holonomy (over the
Novikov field over $\K=\C$). In this situation, we define
$$CF((L_0,\mathcal{E}_0),(L_1,\mathcal{E}_1))=\bigoplus_{p\in
\mathcal{X}(L_0,L_1)} \hom(\mathcal{E}_{0|p},\mathcal{E}_{1|p}),$$
and modify the definition of $\mu^k$ as follows.
Fix objects $(L_0,\mathcal{E}_0),\dots,(L_k,\mathcal{E}_k)$, 
intersections $p_1,\dots,p_k,q$, and a homotopy class
$[u]$. Set $p_0=p_{k+1}=q$ for simplicity. 
Parallel transport along the portion of the boundary of 
$[u]$ that lies on $L_i$ yields an isomorphism
$\gamma_i\in\hom(\mathcal{E}_{i|p_i},\mathcal{E}_{i|p_{i+1}})$ for each
$i=0,\dots,k$.
Now, given elements $\rho_i\in
\hom(\mathcal{E}_{i-1|p_i},\mathcal{E}_{i|p_i})$ $(i=1,\dots,k)$, the composition of
all these linear maps defines an element 
$\eta_{[u],\rho_k,\dots,\rho_1}=\gamma_k\cdot \rho_k\cdot \dots\cdot \gamma_1\cdot \rho_1 \cdot \gamma_0
\in \hom(\mathcal{E}_{0|q},\mathcal{E}_{k|q})$.
Then we set
$$\mu^k(\rho_k,\dots,\rho_1)=\sum_{\substack{q\in \mathcal{X}(L_0,L_k)\\
[u]:\mathrm{ind}([u])=2-k}}
(\#\mathcal{M}(p_1,\dots,p_k,q;[u],J))
T^{\omega([u])}\,\eta_{[u],\rho_k,\dots,\rho_1}.
$$
\end{remark}

\begin{remark}
It is in principle possible to lift the assumption $[\omega]\cdot
\pi_2(M,L)=0$ we have made throughout, at the expense of considerable
analytic and algebraic difficulties in situations where disc bubbling
occurs. Analytically, disc bubbles pose transversality problems that cannot
be solved with the techniques we have described above. Algebraically, they
lead to a {\em curved} $A_\infty$-category, i.e.\ for each object $L$ we
have an element $\mu^0_L\in \hom(L,L)$ which encodes a weighted count of
$J$-holomorphic discs bounded by $L$. The $A_\infty$-relations
\eqref{eq:ainfty} are then modified by allowing the case $\ell=0$ in the
sum. For example, the relation for $k=1$ becomes
$$\mu^1(\mu^1(p))+(-1)^{\deg p}\mu^2(\mu^0_{L_1},p)+\mu^2(p,\mu^0_{L_0})=0,$$
where the last two terms correspond to disc bubbling along either edge of
an index 2 strip. To regain some sanity, one usually considers not arbitrary
objects, but {\em weakly unobstructed} objects, i.e.\ those for which
$\mu^0_L$ is a scalar multiple of the (cohomological) unit of $\hom(L,L)$
(this multiple is sometimes called ``central charge'' or ``superpotential'' 
in the context of mirror symmetry); this happens for instance when
the minimal Maslov index of a holomorphic disc with boundary on $L$ is equal
to two and
Maslov index 2 discs are regular. Weakly unobstructed objects of fixed
central charge then form an honest $A_\infty$-category. The curious reader is
referred to \cite{FO3book}.
\end{remark}

\section{Exact triangles and generators}

While it is usually
impossible to classify all Lagrangian submanifolds of a given symplectic
manifold, or even to directly compute Floer cohomology for all those we can find,
it is often possible to understand the whole Fukaya category in terms of
a small subset of {\em generating objects} -- {\em provided that we understand 
not only differentials and products but also higher operations among those 
generators}. To understand how this comes about, a healthy dose of
homological algebra is necessary; in this section we give a very brief and
informal overview of exact triangles, twisted complexes and generators, in
general and as they pertain to Fukaya categories in particular.
The first part of \cite{SeBook} fills in the many details that we omit here,
and more.

\subsection{Exact triangles and mapping cones}

An exact triangle \bigskip
\begin{equation*}
\begin{array}{c@{\hspace{2cm}}c}
\Rnode{N1}{A} &
\Rnode{N2}{B} \\[1.2cm]
&\Rnode{N3}{C}
\end{array}
\psset{nodesep=0.15cm}
\ncLine{->}{N1}{N2}\Aput{f}
\ncLine{->}{N3}{N1}\Aput{h}\Bput{\smash{\raisebox{-3pt}{$\scriptstyle [1]$}}}
\ncline{->}{N2}{N3}\Aput{g}
\end{equation*}
in an $A_\infty$-category $\mathcal{A}$ consists of a triple of objects $A,B,C$ and
closed morphisms $f\in \hom^0(A,B)$, $g\in \hom^0(B,C)$, $h\in \hom^1(C,A)$
such that $C$ is (up to quasi-isomorphism) a {\em mapping cone} of $f:A\to
B$, with $g$ and $h$ the natural maps to and from it. 
We will clarify the meaning of this definition in the next section; for
now, we simply mention some key features and motivate the concept. 

Exactness means that the
compositions $\mu^2(g,f)$, $\mu^2(h,g)$ and $\mu^2(f,h)$ are exact, i.e.\
in the cohomology category $H(\mathcal{A})$ the maps compose to zero.
(However, their triple Massey product is typically nontrivial.)
An exact triangle induces long exact sequences
on morphism spaces in the cohomology category: for every test object $T$, we have a long exact sequence
\begin{equation}\label{eq:leshom}
\dots\to H^i\hom(T,A)\stackrel{f}{\longrightarrow}
H^i\hom(T,B)\stackrel{g}{\longrightarrow}H^i\hom(T,C)
\stackrel{h}{\longrightarrow}H^{i+1}\hom(T,A)\stackrel{f}{\to}\dots
\end{equation}
where $H^i\hom(T,A)$ is the cohomology of $\hom(T,A)$ with respect to the
differential~$\mu^1$,
and the maps are given by composition (in 
the cohomology category) with $f,g,$ and $h$; and similarly (in the
contravariant direction) for morphisms from $A,B,C$ to~$T$.
Moreover, as $T$ varies these long exact sequences fit together
naturally with respect to the multiplicative action of the groups 
$H^*\hom(T',T)$, i.e.\ \eqref{eq:leshom} fits into an exact sequence of
modules over $H(\mathcal{A})$.

Exact triangles can also be characterized as images under
$A_\infty$-functors of a ``universal'' abstract exact triangle 
living in an $A_\infty$-category with three objects~\cite[\S 3g]{SeBook}.

The $A_\infty$-category $\mathcal{A}$ is said to be {\em triangulated}\/ if
every closed morphism $f:A\to B$ can be completed to an exact triangle
(and the shift functor [1] acting on $\mathcal{A}$ by change of
gradings is a quasi-equivalence); or, in other terms, if all morphisms in
$\mathcal{A}$ have mapping cones. Here it is important to point out a key
difference with the case of ordinary triangulated categories, where the
triangles are an additional piece of structure on the category: 
the $A_\infty$-structure is rich enough to ``know'' about triangles,
and triangles automatically satisfy an analogue of the usual axioms.
In the same vein, $A_\infty$-functors are always exact, i.e.\ map exact
triangles to exact triangles.

Before saying more about mapping cones in $A_\infty$-categories, let us
discuss some classical motivating examples.

\begin{example}
The mapping cone of a continuous map $f:X\to Y$ between topological spaces
is, by definition, the space obtained from $X\times [0,1]$ by attaching
$Y$ to $X\times \{1\}$ via the map $f$ and collapsing $X\times\{0\}$ to a
point:
$$\mathrm{Cone}(f)=\bigl((X\times [0,1]) \sqcup Y\bigr) \big/ (x,0)\sim (x',0),\ (x,1)\sim
f(x)\ \forall x,x'\in X.$$
We then have a sequence of maps
$$X\stackrel{f}{\longrightarrow} Y\stackrel{i}{\longrightarrow}
\mathrm{Cone}(f)\stackrel{p}{\longrightarrow} \Sigma X\stackrel{\Sigma
f}{\longrightarrow}\Sigma Y\to\dots,$$
where $i$ is the inclusion of $Y$ into the mapping cone, and $p$ is the
projection to the suspension of $X$ obtained by collapsing $Y$. 
The composition of any two of these maps is nullhomotopic, and the induced
maps on (co)homology form a long exact sequence.
\end{example}

\begin{example}
The notion of mapping cone in the category of chain complexes is directly
modelled on the previous example: let $A=(\bigoplus A^i, d_A)$ and 
$B=(\bigoplus B^i, d_B)$ be two chain complexes, and let $f:A\to B$ be
a chain map (i.e., a collection of maps $f^i:A^i\to B^i$ satisfying
$d_B f^i+ f^{i+1} d_A=0$). Then the mapping cone of $f$ is, by definition,
the chain complex $C=A[1]\oplus B$ (i.e., $C^i=A^{i+1}\oplus
B^i$), equipped with the differential $$d_C=\begin{pmatrix} d_A& 0\\
f&d_B\end{pmatrix}.$$
The map $f$, the inclusion of $B$ into $C$ as a subcomplex, and the
projection of $C$ onto the quotient complex $A[1]$ then fit into an exact
sequence.
\end{example}

\begin{example}\label{ex:modulecone}
Let $A$ be an algebra (resp.\ differential graded algebra or
$A_\infty$-algebra), and consider the category of differential graded
modules (resp.\ $A_\infty$-modules) over~$A$. Recall that such a module $M$ 
is a chain complex equipped with a degree 1 differential $d_M$
and a multiplication map $A\otimes M\to M$,
$(a,m)\mapsto a\cdot m$, satisfying the Leibniz rule and associative
(up to homotopies given by higher structure maps \hbox{$\mu^{\smash{k|1}}_M:A^{\otimes k}
\otimes M\to M[1-k]$,} in the case of $A_\infty$-modules).
The mapping cone of a module homomorphism $f:M\to N$ can then be defined
essentially as in the previous example. In the differential graded case,
$f$ is a chain map compatible with the multiplication, and the mapping
cone of $f$ as a chain complex inherits a natural module structure.
For $A_\infty$-modules, recalling that
an $A_\infty$-homomorphism is a collection of maps $f^{k|1}:A^{\otimes
k}\otimes M\to N[-k]$ (where the linear term $f^{0|1}$ is a chain map compatible
with the product $\mu^{1|1}$ up to a homotopy given by $f^{1|1}$, and so on),
the structure maps $\mu^{k|1}_K:A^{\otimes k}\otimes K\to K[1-k]$
($k\ge 0$) of the mapping cone $K=M[1]\oplus N$ are 
given by
$$\mu^{k|1}_K(a_1,\dots,a_k,(m,n))=(\mu^{k|1}_M(a_1,\dots,a_k,m),
f^{k|1}(a_1,\dots,a_k,m)+\mu^{k|1}_N(a_1,\dots,a_k,n)).$$
\end{example}

\subsection{Twisted complexes}
When an $A_\infty$-category $\mathcal{A}$ is not known to be triangulated, 
it is often advantageous to embed it into a larger category in which mapping
cones are guaranteed to exist. For example, one can always do so by
using the {\em Yoneda embedding} construction into the category of 
$A_\infty$-modules over $\mathcal{A}$ (in which mapping cones always
exist, cf.\ Example \ref{ex:modulecone}); see e.g.\ \cite[\S 1]{SeBook}.
A milder construction, which retains more features of the original category
$\mathcal{A}$, involves {\em twisted complexes}.
We give a brief outline, and refer the reader to \cite[\S 3]{SeBook} for details.

\begin{definition}\label{def:tw}
A {\em twisted complex} $(E,\delta^E)$ consists of:
\begin{itemize}
\item a formal direct sum $E=\bigoplus\limits_{i=1}^N E_i[k_i]$
of shifted objects of $\mathcal{A}$ $($i.e., a finite collection of pairs
$(E_i,k_i)$ where $E_i\in \mathrm{ob}\,\mathcal{A}$ and $k_i\in \Z);$
\item a strictly lower triangular differential $\delta^E\in \mathrm{End}^1(E)$,
i.e.\ a collection of maps $\delta^E_{ij}\in \mathrm{Hom}^{k_j-k_i+1}(E_i,E_j)$, 
$1\le i<j\le N$,
satisfying the equation 
\begin{equation}\label{eq:maurercartan}
\sum_{k\ge 1} \mu^k(\delta^E,\dots,\delta^E)=0,
\end{equation}
i.e., \ $\sum\limits_{k\ge 1} \sum\limits_{\ i=i_0<i_1<\dots<i_k=j} 
\mu^k(\delta^E_{i_{k-1}i_k},\dots,\delta^E_{i_0i_1})=0$ for all\/ $1\le i<j\le N$.
\end{itemize}

\noindent A degree $d$ morphism of twisted complexes is simply a degree $d$
map between the underlying formal direct sums, i.e.\ if
$E=\bigoplus E_i[k_i]$ and $E'=\bigoplus E'_j[k'_j]$ then an element of
$\mathrm{Hom}^d(E,E')$ is by definition a 
collection of morphisms $a_{ij}\in \mathrm{Hom}^{d+k'_j-k_i}(E_i,E'_j)$.

Finally, given twisted complexes $(E_0,\delta^{0}),\dots,
(E_k,\delta^{k})$,  $k\ge 1$,
and morphisms $a_i\in \mathrm{Hom}(E_{i-1},E_i)$, we set
\begin{equation*}
\mu^k_{\mathrm{Tw}}(a_k,\dots,a_1)=
\sum_{j_0,\dots,j_k\ge 0}
\mu^{k+j_0+\dots+j_k}(\underbrace{\delta^{k},\dots,\delta^{k}}_{j_k},
a_k,\dots,\underbrace{\delta^{1},\dots,\delta^{1}}_{j_1},a_1,
\underbrace{\delta^{0},\dots,\delta^{0}}_{j_0}).
\end{equation*}
(The sum is finite since each $\delta^i$ is strictly lower triangular).
\end{definition}

\begin{proposition}
The above construction defines a triangulated $A_\infty$-category
which we denote by $\mathrm{Tw}\,\mathcal{A}$, and into which $\mathcal{A}$
embeds fully faithfully.
\end{proposition}

It is instructive to see how twisted complexes relate to ordinary chain
complexes:

\begin{example}
Given objects $A,B,C$ of $\mathcal{A}$ and $f\in \hom^0(A,B)$, $g\in
\hom^0(A,C)$, we can consider $(A[2]\oplus B[1]\oplus C,\delta=f+g)$,
conventionally denoted by $$\{A\stackrel{f}{\longrightarrow}
B\stackrel{g}{\longrightarrow} C\}.$$
This forms a twisted complex if and only if $\mu^1(f)=\mu^1(g)=0$ and
$\mu^2(g,f)=0$, i.e.\ $f$ and $g$ are closed morphisms
and their composition is zero. However, we can also introduce
an extra term $h\in \hom^{-1}(A,C)$ into the differential $\delta$, in 
which case the last condition becomes $\mu^2(g,f)+\mu^1(h)=0$: thus it is
sufficient for the composition of $f$ and $g$ to be exact, with a homotopy
given by $h$. 
\end{example}

\begin{definition}
Given twisted complexes $(E,\delta),(E',\delta')\in \mathrm{Tw}\,\mathcal{A}$ 
and a closed morphism $f\in \hom^0(E,E')$ $($i.e., such that
$\mu^1_{\mathrm{Tw}}(f)=0)$, the {\em abstract mapping cone} of $f$ is
the twisted complex $$\mathrm{Cone}(f)=\left(E[1]\oplus E',\begin{pmatrix}
\delta&0\\f&\delta'\end{pmatrix}\right).$$
Given objects $A,B,C$ of $\mathcal{A}$ and a closed morphism $f\in
\hom^0(A,B)$, we say that $C$ is a {\em mapping cone} of $f$ if, in the
category of twisted complexes $\mathrm{Tw}\,\mathcal{A}$, the object $C$
is quasi-isomorphic to the abstract mapping cone of $f$,
$\{A\stackrel{f}{\longrightarrow}B\}=(A[1]\oplus B,f)$.
\end{definition}

When $C$ is a mapping cone of $f:A\to B$,
by composing the inclusion of $B$ into
the abstract mapping cone (resp.\ the projection to $A[1]$) with the given
quasi-isomorphism from the abstract mapping cone to $C$ (resp.\ its
quasi-inverse) we obtain morphisms $i:B\to C$ and $p:C\to A[1]$, which
sit with $f$ in an exact triangle.

\subsection{Exact triangles in the Fukaya category}
The reader may legitimately wonder about the relevance of the above 
discussion to Fukaya categories. It turns out that at least some
mapping cones in the Fukaya category of a symplectic manifold can be
understood geometrically. There are two well-known sources of these:
Dehn twists, and Lagrangian connected sums.

\subsubsection{Dehn twists} 
\begin{figure}[t]
\setlength{\unitlength}{5mm}
\begin{picture}(7,5)(0,-1)
\psset{unit=\unitlength}
\psline{->}(0,0)(6,0)
\psline{->}(0,0)(0,4)
\put(5.5,0.5){\small $\|p\|$}
\put(0.3,3.7){\small $h$}
\psline[linestyle=dotted](0,3)(4,3)
\pscurve(0,0)(0.7,1)(1.45,2)(2.6,2.9)(4,3)(5.8,3)
\psline[linewidth=0.5pt](0.5,0)(0.5,0.714)
\put(0.6,0.2){\tiny $\pi$}
\end{picture}
\qquad\qquad
\begin{picture}(12,6)(-2,-3)
\psset{unit=\unitlength}
\psellipticarc(0,-2.5)(2,0.5){-180}{0}
\psellipticarc[linestyle=dashed](0,-2.5)(2,0.5){0}{180}
\psellipse(0,2.5)(2,0.5)
\psline(-2,-2.5)(-2,2.5)
\psline(2,-2.5)(2,2.5)
\psellipticarc(0,0)(2,0.5){-180}{0}
\psellipticarc[linestyle=dashed](0,0)(2,0.5){0}{180}
\put(-2.4,0){\makebox(0,0)[rc]{$S$}}
\psline(0,-3)(0,2)
\pscircle*(0,-0.5){0.15}
\psline{|->}(3,0)(5,0)
\put(4,0.2){\makebox(0,0)[cb]{$\tau_S$}}
\psellipticarc(8,-2.5)(2,0.5){-180}{0}
\psellipticarc[linestyle=dashed](8,-2.5)(2,0.5){0}{180}
\psellipse(8,2.5)(2,0.5)
\psline(6,-2.5)(6,2.5)
\psline(10,-2.5)(10,2.5)
\psellipticarc(8,0)(2,0.5){-180}{0}
\psellipticarc[linestyle=dashed](8,0)(2,0.5){0}{180}
\pscurve(8,-3)(8.1,-2)(9.9,-0.8)(10,-0.5)
\pscurve[linestyle=dashed](10,-0.5)(9.9,-0.2)(8,0.5)(6.5,0.7)(6,0.6)(6,0.5)
\pscurve(6,0.5)(6.5,0.5)(7.9,1.3)(8,2)
\pscircle*(8,0.5){0.15}
\end{picture}
\caption{The generating Hamiltonian on the complement of the zero section
in $T^*S^n$, and the action of the Dehn twist on a cotangent fiber.}
\label{fig:dehntwist}
\end{figure}
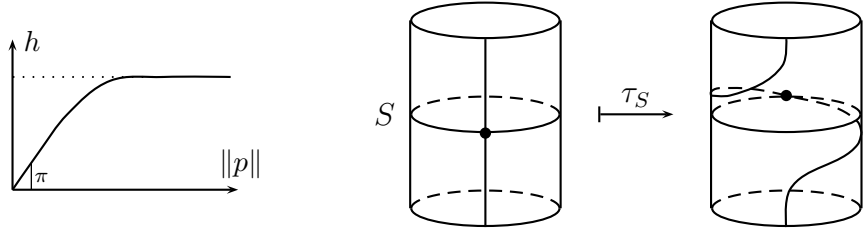

The symplectic geometry of Dehn twists was first
considered by Arnold,
and later studied extensively by Seidel \cite{SeLES,SeBook}.
The local model is as follows. In the cotangent bundle $T^*S^n$ equipped
with its canonical symplectic form, a Hamiltonian of the form
$H(p,q)=h(\|p\|)$ (where $p$ is the fiber coordinate and $\|\cdot\|$ is the
standard metric) generates a rescaled version of geodesic flow. Choosing 
$h:[0,\infty)\to\R$ so that $h'(0)=\pi$, $h''\le 0$, and $h$ is constant 
outside of a neighborhood of zero, we obtain a Hamiltonian diffeomorphism
of the complement of the zero section $T^*S^n\setminus S^n$, which
can be extended across the zero section by the antipodal 
map on $S^n$ to obtain a symplectomorphism of $T^*S^n$ (see Figure \ref{fig:dehntwist}).

Now, given a Lagrangian sphere $S$ in a symplectic manifold $(M,\omega)$,
by Weinstein's theorem a neighborhood of $S$ in $M$ is symplectomorphic
to a neighborhood of the zero section in $T^*S^n$; thus, performing the
above construction inside the standard neighborhood of $S$, we obtain a
symplectomorphism $\tau_S$, the {\em Dehn twist} about $S$, which is
supported in a neighborhood of $S$ and maps $S$ to itself antipodally. 
(Note: $\tau_S$ depends on the choices
made in the construction, but its isotopy class doesn't.)

\begin{theorem}[Seidel \cite{SeLES,SeBook}]
Given a Lagrangian sphere $S$ and any object $L$ of $\F(M,\omega)$,
there is an exact triangle in $\mathrm{Tw}\,\F(M,\omega)$,\bigskip
\begin{equation}\label{eq:striangle}
\begin{array}{c@{\hspace{1cm}}c}
\Rnode{N1}{HF^*(S,L)\otimes S} &
\Rnode{N2}{L} \\[1.2cm]
&\Rnode{N3}{\tau_S(L)}
\end{array}
\psset{nodesep=0.15cm}
\everypsbox{\scriptstyle}
\ncLine{->}{N1}{N2}\Aput{\qquad\quad ev}
\ncLine{->}{N3}{N1}\Aput{[1]}
\ncline{->}{N2}{N3}
\end{equation}
In other terms, the object
$\tau_S(L)$ of $\F(M,\omega)$ is quasi-isomorphic in 
$\mathrm{Tw}\,\F(M,\omega)$ to the abstract mapping cone of $ev$.
\end{theorem}
\noindent
In \eqref{eq:striangle}, $HF^*(S,L)\otimes S$ is a direct sum of shifted copies of $S$,
with one summand for each generator of $HF^*(S,L)$, and $ev$ is a
tautological evaluation map, mapping each summand to $L$ by a closed morphism 
representing the given generator of $HF^*(S,L)=H^*\mathrm{Hom}(S,L)$.

Given a test object $T$, the corresponding long exact sequence \eqref{eq:leshom}
is Seidel's long exact sequence in Floer cohomology \cite{SeLES} associated to the Dehn
twist $\tau_S$ for all $T,L$:
\begin{equation*}
\dots\to HF^*(S,L)\otimes HF^*(T,S)\stackrel{\mu^2}{\longrightarrow}
HF^*(T,L)\longrightarrow
HF^*(T,\tau_S(L))\stackrel{[1]}{\longrightarrow}\dots
\end{equation*}

\subsubsection{Lagrangian connected sums}
Given two Lagrangian submanifolds $L_1,L_2$ which intersect
transversely in a single point $p$,
we can form the {\em Lagrangian connected sum} (or {\em surgery}
in the terminology of \cite{Pol} and \cite{FO3surgery}) $L_1\# L_2$.
One possible construction is as follows. For $\epsilon>0$, the graph of 
the 1-form $\epsilon\, d\log \|x\|$ on $\R^n$,
given by the equations $y_i=\epsilon\,x_i/\|x\|^2$, 
is a Lagrangian submanifold of $T^*\R^n\simeq \C^n$ which is asymptotic to 
the zero section (i.e., $\R^n\subset\C^n$) as $\|x\|\to \infty$ and to the cotangent fiber 
over zero (i.e., $(i\R)^n\subset \C^n$) as $\|y\|\to \infty$; using suitable cut-off functions, we can modify
this Lagrangian so that it agrees with $\R^n\cup (i\R)^n$
outside of a small neighborhood of the origin. Pasting this local model into
a suitable Darboux chart centered at the intersection point $p$ and chosen
so that $T_pL_1=\R^n$ and $T_pL_2=(i\R)^n$ yields
$L_1\# L_2$.
(Note that, for a single connected sum operation, the end result is
independent of the size parameter $\epsilon$ and other choices up to 
Hamiltonian isotopy; not so when summing at multiple points. Also note that
$L_2\# L_1$ is {\em not} isotopic to $L_1\# L_2$.)

\begin{remark}
When $L_2$ is a sphere, $L_1\# L_2$ is Hamiltonian isotopic to
$\tau_{L_2}(L_1)$; this provides the basis for an alternative description
of the connected sum operation.
\end{remark}

\begin{figure}[t]
\setlength{\unitlength}{20mm}
\begin{picture}(2,2)(-0.8,-0.8)
\psset{unit=\unitlength}
\newgray{gray15}{0.9}
\psline[linearc=0.4](-0.8,0)(0,0)(0,-0.8)
\psline[linearc=0.4,linestyle=none,fillcolor=gray15,fillstyle=solid](0.9,0)(0,0)(0,0.9)
\psline[linearc=0.4](1.2,0)(0,0)(0,1.2)
\psline(-0.3,1.2)(1.2,-0.3)
\put(-0.65,0.4){\small $L_1\# L_2$}
\put(0.52,0.5){\small $T$}
\end{picture}
\qquad\qquad
\begin{picture}(2,2)(-0.8,-0.8)
\psset{unit=\unitlength}
\newgray{gray15}{0.9}
\pspolygon[fillstyle=solid,fillcolor=gray15](0,0)(0,0.9)(0.9,0)
\psline(-0.8,0)(1.2,0)
\psline(0,-0.8)(0,1.2)
\psline(-0.3,1.2)(1.2,-0.3)
\pscircle*(0,0){0.035}
\put(-0.15,-0.15){\small $p$}
\put(0.35,-0.2){\small $L_1$}
\put(-0.28,0.4){\small $L_2$}
\put(0.52,0.5){\small $T$}
\end{picture}
\caption{The Lagrangian connected sum $L_1\# L_2$ vs.\ $L_1\cup L_2$}
\label{fig:surgery}
\end{figure}
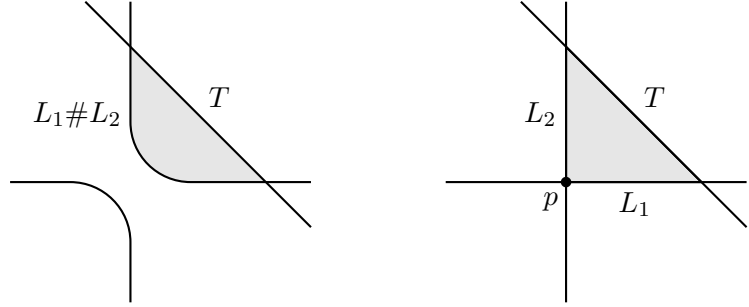

Given some other Lagrangian submanifold $T$ (in generic position relatively
to $L_1$ and $L_2$), choosing $\epsilon$ small enough in the above
construction ensures that the intersections of $T$ with $L_1\# L_2$
are the same as with $L_1\cup L_2$. 
Fukaya-Oh-Ohta-Ono \cite{FO3surgery} have studied the moduli spaces
of $J$-holomorphic discs bounded by $L_1\# L_2$ and $T$.
Their main result is that, for suitable $J$ and small 
enough $\epsilon$, $J$-holomorphic strips with boundary on $T$ and $L_1\#
L_2$ connecting an intersection in $T\cap L_2$ to one in $T\cap L_1$ are
in bijection with $J$-holomorphic triangles bounded by $T$, $L_2$ and $L_1$
with a corner at $p$, whereas the counts of rigid strips in the other
direction vanish. This is elementary in dimension 1, as illustrated
by Figure \ref{fig:surgery}, but much harder in higher dimensions. 

The outcome is that, as a chain complex, $CF(T,L_1\# L_2)$ is the mapping cone
of the map $\mu^2(p,\cdot):CF(T,L_2)\to CF(T,L_1)$ given by multiplication by
the generator $p$ of $CF(L_2,L_1)$. Hence, the short exact sequence
$$0\to CF(T,L_1) \to CF(T,L_1\# L_2)\to CF(T,L_2)\to 0$$
induces a long exact sequence 
$$\dots\to HF(T,L_1)\longrightarrow HF(T,L_1\# L_2)
\longrightarrow HF(T,L_2)\xrightarrow{\mu^2([p],\cdot)}
HF(T,L_1)\to \dots$$

By an analogous argument for higher structure maps, one expects that
this long exact sequence can be upgraded to an exact triangle in the
Fukaya category,\medskip
\begin{equation}\label{eq:surgtriangle}
\begin{array}{c@{\hspace{1cm}}c}
\Rnode{N1}{L_2} &
\Rnode{N2}{L_1} \\[1cm]
&\Rnode{N3}{L_1\# L_2}
\end{array}
\psset{nodesep=0.15cm}
\everypsbox{\scriptstyle}
\ncLine{->}{N1}{N2}\Aput{p}
\ncLine{->}{N3}{N1}\Aput{[1]}
\ncline{->}{N2}{N3}
\end{equation}
i.e., $L_1\# L_2$ is quasi-isomorphic to the twisted complex
$\mathrm{Cone}(p)=\{L_2\overset{p}{\to}L_1\}$.

\noindent (If $L_2$ is a sphere, this is Seidel's exact triangle for the Dehn
twist of $L_1$ about $L_2$.)

\begin{remark}
Recall that, by definition, the differential
$\mu^1_{\mathrm{Tw}}$ on $\hom(T,\mathrm{Cone}(p))$ involves not
only the original Floer differential $\mu^1$, but also multiplication
by the differential of the twisted complex, i.e.\ $\mu^2(p,\cdot)$.
This is exactly consistent with the above description
of $J$-holomorphic strips with boundary on $T$ and $L_1\# L_2$.
Thus, replacing Lagrangian submanifolds by quasi-isomorphic 
twisted complexes built out of simpler Lagrangians, while computationally 
powerful, comes at the expense of having to consider higher operations
on their Floer complexes (in this case, the expression for 
$\mu^1_{\mathrm{Tw}}$ involves $\mu^2$, and similarly that for 
$\mu^2_{\mathrm{Tw}}$ involves $\mu^3$).
\end{remark}

\subsection{Generation and Yoneda embedding}
\subsubsection{Generators and split-generators}
\begin{definition}
The objects $G_1,\dots,G_r$ are said to {\em generate} the $A_\infty$-category
$\mathcal{A}$ if, in $\mathrm{Tw}\,\mathcal{A}$, every object of
$\mathcal{A}$ is quasi-isomorphic to a twisted complex built from copies
of $G_1,\dots,G_r$. $($In other terms, every object of $\mathcal{A}$ can be
obtained from $G_1,\dots,G_r$ by taking iterated mapping cones.$)$

The objects $G_1,\dots,G_r$ are said to {\em split-generate} $\mathcal{A}$ if every
object of $\mathcal{A}$ is quasi-isomorphic to a direct summand in a
twisted complex built from copies of $G_1,\dots,G_r$.
\end{definition}
\begin{figure}[b]
\setlength{\unitlength}{35mm}
\begin{picture}(1,1)(0,0)
\psset{unit=\unitlength}
\newgray{gray15}{0.9}
\psline(0,0)(1,0)(1,1)(0,1)(0,0)
\psline{->>}(0,0)(0,0.45)
\psline{->>}(1,0)(1,0.45)
\psline{->}(0,1)(0.45,1)
\psline{->}(0,0)(0.45,0)
\psline(0.55,0)(0.55,1)
\psline(0,0.55)(1,0.55)
\psline[linearc=0.3](0,0.52)(0.45,0.6)(0.52,1)
\psline[linearc=0.3](0.52,0)(0.6,0.45)(1,0.52)
\put(0.58,0.8){\small$\beta$}
\put(0.75,0.58){\small$\alpha$}
\pscircle*(0.55,0.55){0.02}
\put(0.47,0.47){\small $p$}
\end{picture}
\qquad\qquad
\begin{picture}(1,1)(0,0)
\psset{unit=\unitlength}
\newgray{gray15}{0.9}
\psline(0,0)(1,0)(1,1)(0,1)(0,0)
\psline{->>}(0,0)(0,0.45)
\psline{->>}(1,0)(1,0.45)
\psline{->}(0,1)(0.45,1)
\psline{->}(0,0)(0.45,0)
\psline(0.5,0)(0.5,1)
\psline(0,1)(1,0.5)
\psline(0,0.5)(1,0)
\psline[linearc=0.05](0,1)(0.47,0.78)(0.5,1)
\psline[linearc=0.05](0.5,0.5)(0.53,0.72)(1,0.5)
\psline[linearc=0.15](0,0.5)(0.47,0.28)(0.5,0.5)
\psline[linearc=0.15](0.5,0)(0.53,0.22)(1,0)
\put(0.53,0.35){\small$\beta$}
\put(0.8,0.65){\small$\gamma$}
\pscircle*(0.5,0.75){0.02}
\put(0.39,0.67){\small $q_1$}
\pscircle*(0.5,0.25){0.02}
\put(0.39,0.17){\small $q_2$}
\end{picture}
\caption{Split-generating the Fukaya category of $T^2$}
\label{fig:genT2}
\end{figure}
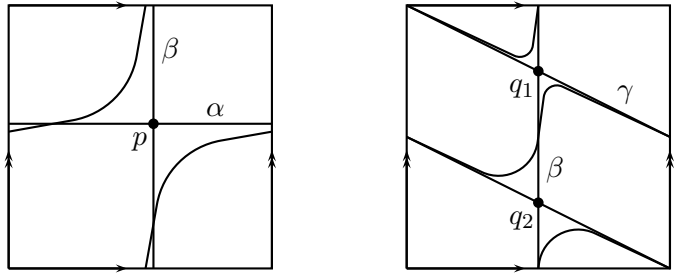
\begin{example}
Consider the Fukaya
category of the torus $T^2$ with its standard area form.
Starting from the standard curves $\alpha$ and $\beta$ along the two
factors of the torus, by taking iterated mapping cones we can obtain
simple closed curves representing all nontrivial primitive elements in
$\pi_1(T^2)=\Z^2$.
For instance, the loop $\tau_\alpha(\beta)\simeq \beta\#\alpha$
(Figure \ref{fig:genT2} left) is quasi-isomorphic to the
mapping cone of $p\in \mathrm{Hom}(\alpha,\beta)$;
further applications of the Dehn twists $\tau_\alpha$ and $\tau_\beta$ 
(which generate the mapping class group of $T^2$) eventually yield
simple closed curves in all primitive homotopy classes. However, the
objects obtained in this manner all satisfy a certain ``balancing''
condition: given a 1-form $\theta\in \Omega^1(T^2\setminus \{pt\})$ with
$d\theta=\omega$ and such that $\int_\alpha\theta=\int_\beta\theta=0$,
$\theta$ also integrates to zero on all iterated mapping cones built from
$\alpha$ and $\beta$. For instance, all the simple closed curves that can be 
obtained in a given homotopy class are Hamiltonian isotopic to each 
other. Thus, $\alpha$ and $\beta$ generate the subcategory of
$\F(T^2)$ consisting of Lagrangians which are balanced with respect to
$\theta$, but not all of $\F(T^2)$. 

On the other hand, given the two loops $\beta$ and $\gamma$ shown
on Figure \ref{fig:genT2} right, the mapping cone of $T^{a_1} q_1+T^{a_2}
q_2\in \mathrm{Hom}(\gamma,\beta)$ can be interpreted geometrically as the
connected sum of $\beta$ and $\gamma$ at their two intersection points
$q_1$ and $q_2$, with different gluing parameters.
This mapping cone is therefore quasi-isomorphic to the direct sum
of two simple closed curves in the homotopy class of $\alpha$, but whose
Hamiltonian isotopy classes depend on $a_1$ and $a_2$. Thus, by considering
direct summands in mapping cones we can obtain all nontrivial simple closed 
curves up to Hamiltonian isotopy, rather than only those that are balanced:
$\alpha$ and $\beta$ split-generate $\F(T^2)$.
\end{example}

\subsubsection{Yoneda embedding}\label{ss:yoneda}
Let $G_1,\dots,G_r$ be split-generators of the $A_\infty$-category
$\mathcal{A}$. Then the endomorphism algebra of $G_1\oplus\dots\oplus G_r$,
$$\mathcal{G}=\bigoplus_{i,j=1}^r \hom(G_i,G_j)$$
is an $A_\infty$-algebra (with structure maps given by the operations
$\mu^k$ of $\mathcal{A}$). Next, given any object $L$ of $\mathcal{A}$,
$$\mathcal{Y}(L)=\bigoplus_{i=1}^r \hom(G_i,L)$$
is a (right) $A_\infty$-module over $\mathcal{G}$, with differential
given by $\mu^1$, multiplication $\mu^{1|1}$ given by the operations
$$\hom(G_j,L)\otimes \hom(G_i,G_j)\stackrel{\mu^2}{\longrightarrow}
\hom(G_i,L),$$ and so on (the structure map $\mu^{1|k}$ of $\mathcal{Y}(L)$
is given by $\mu^{k+1}$).

Moreover, to a morphism $a\in \hom(L,L')$ we can associate an
$A_\infty$-homomorphism $\mathcal{Y}(a)\in \hom_{\text{mod-}\mathcal{G}}(\mathcal{Y}(L),
\mathcal{Y}(L'))$, whose linear term is given by composition with $a$.

The assignment $L\mapsto \mathcal{Y}(L)$, $a\mapsto \mathcal{Y}(a)$ is in turn the linear term of
an $A_\infty$-functor $\mathcal{Y}$, which is the restriction to the
given set of objects $G_1,\dots,G_r$ of the $A_\infty$ {\em Yoneda embedding} 
$\mathcal{A}\to \text{mod-}\mathcal{A}$ (see e.g.\ \cite[\S 1]{SeBook}):

\begin{proposition}
The above construction extends to an $A_\infty$-functor
$\mathcal{Y}$ from $\mathcal{A}$ to $\text{mod-}\mathcal{G}$.
Moreover, if $G_1,\dots,G_r$ split-generate $\mathcal{A}$ then
this  $A_\infty$-functor is a fully faithful quasi-embedding.
\end{proposition}

\section{The wrapped Fukaya category, examples and applications}

In this section we assume that $(M,\omega)$ is a {\em Liouville manifold},
i.e.\ an exact symplectic manifold such that the Liouville vector field $Z$
associated to the chosen primitive $\theta\in \Omega^1(M)$ of the symplectic 
form (i.e., the conformally symplectic vector field defined by 
$\iota_Z\omega=\theta$) is complete and outward pointing at infinity.
More precisely, we require that $M$ contains a compact domain $M^{in}$ with 
boundary a smooth hypersurface $\partial M$ on which
$\alpha=\theta_{|\partial M}$ is a contact form, and $Z$ is positively transverse to
$\partial M$ and has no zeroes outside of $M^{in}$. The flow of $Z$ can then
be used to identify $M\setminus M^{in}$ with
the positive symplectization $(1,\infty)\times \partial M$
equipped with the exact symplectic form $\omega=d(r\alpha)$
and the Liouville field $Z=r\frac{\partial}{\partial r}$.

In this setting it is natural to consider not only compact exact Lagrangian
submanifolds as we have done above, but also some noncompact ones with
suitable behavior at infinity. There are two different types of such 
noncompact Fukaya categories, depending on the manner in which perturbations
at infinity are used to define Floer complexes. One possibility is to perform
``small'' perturbations at infinity, restricting oneself to a smaller set
of ``admissible'' objects which go to infinity along well-controlled 
directions. Two constructions that follow this philosophy are the ``infinitesimal'' Fukaya category first 
defined by Nadler and Zaslow for cotangent bundles \cite{NZ} and later
extended to Liouville manifolds equipped with a choice of Lagrangian
skeleton; and Fukaya categories of Lefschetz fibrations as constructed
by Seidel \cite{SeBook,SeLF}, and their putative generalization
to Landau-Ginzburg models, in which the behavior at infinity is controlled
by a projection to the complex plane. Here we focus on the other approach,
which is to consider {\em large} perturbations at infinity, leading to
the {\em wrapped Fukaya category} of Abouzaid and Seidel
\cite{AS,AbGenerate}. For completeness we mention the nascent subject of
{\em partially wrapped} Fukaya categories, which attempt to interpolate between
these two approaches (cf.\ e.g.\ \cite{Aur10}).

\subsection{The wrapped Fukaya category}

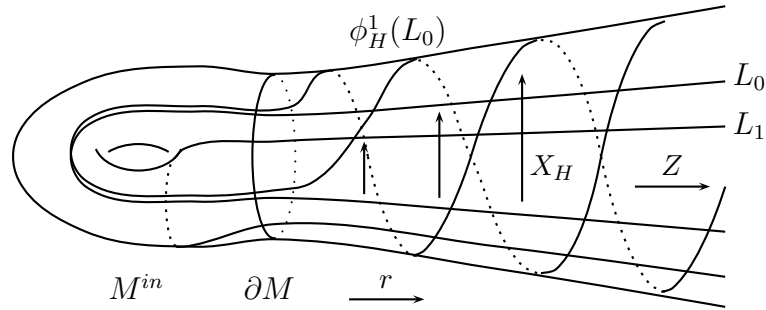
\begin{figure}[b]
\setlength{\unitlength}{1cm}
\begin{picture}(9.5,4)(-1.5,-2)
\psset{unit=\unitlength,dash=1pt 2.5pt}
\pscurve(8,2)(5,1.5)(2,1.1)(1,1.2)(-0.5,1)(-1.4,0.3)(-1.4,-0.3)(-0.5,-1)(1,-1.2)(2,-1.1)(5,-1.5)(8,-2)
\psellipticarc(2,0)(0.3,1.1){90}{270}
\psellipticarc[linestyle=dotted](2,0)(0.3,1.1){-90}{90}
\psellipticarc(0.2,0.2)(0.6,0.4){-170}{-10}
\psellipticarc(0.2,-0.22)(0.6,0.4){35}{145}
\put(-0.2,-1.9){$M^{in}$}
\put(1.6,-1.9){$\partial M$}
\put(3.4,-1.75){\small $r$}
\psline{->}(3,-1.9)(4,-1.9)
\pscurve(8,1)(5,0.75)(2,0.55)(1,0.6)(-0.5,0.4)(-0.7,0)(-0.5,-0.4)(1,-0.6)(2,-0.55)(5,-0.75)(8,-1)
\pscurve(2.8,1.15)(2.6,1.1)(2.3,0.7)(2,0.62)(1,0.67)(-0.5,0.48)(-0.7,0.07)(-0.5,-0.32)(1,-0.52)(2,-0.45)(2.5,-0.35)(3.25,0.5)(3.7,1.2)(3.9,1.3)
\pscurve[linestyle=dashed](2.75,1.15)(2.85,1.1)(3.25,0)
\pscurve[linestyle=dashed](3.25,0)(3.8,-1.3)(3.9,-1.32)
\pscurve(3.85,-1.32)(4,-1.25)(4.62,0)
\pscurve(4.62,0)(5.3,1.45)(5.5,1.55)
\pscurve[linestyle=dashed](5.55,1.55)(5.75,1.5)(6.35,0)
\pscurve[linestyle=dashed](6.35,0)(7,-1.7)(7.2,-1.8)
\pscurve(7.15,-1.8)(7.3,-1.75)(8,-0.4)
\pscurve[linestyle=dashed](3.8,1.3)(4,1.25)(4.62,0)(5.3,-1.45)(5.5,-1.55)
\pscurve(5.55,-1.55)(5.75,-1.5)(6.35,0)
\pscurve(6.35,0)(7,1.7)(7.2,1.8)
\pscurve(8,0.4)(5,0.3)(2,0.2)(0.8,0.1)(0.7,0)
\psellipticarc[linestyle=dashed](0.7,-0.6)(0.15,0.6){90}{270}
\pscurve(0.7,-1.2)(0.8,-1.2)(2,-0.9)(5,-1.2)(8,-1.6)
\put(3,1.55){$\phi_H^1(L_0)$}
\put(8.1,0.95){$L_0$}
\put(8.1,0.3){$L_1$}
\psline{->}(3.2,-0.5)(3.2,0.2)
\psline{->}(4.2,-0.55)(4.2,0.6)
\psline{->}(5.3,-0.6)(5.3,1.1)
\psline{->}(6.8,-0.4)(7.8,-0.4)
\put(7.15,-0.25){\small $Z$}
\put(5.4,-0.25){\small $X_H$}
\end{picture}
\caption{Wrapping by a quadratic Hamiltonian}
\label{fig:wrap}
\end{figure}

The objects of the wrapped Fukaya category $\mathcal{W}(M)$ of a Liouville manifold 
$(M,\omega=d\theta)$ are exact Lagrangian submanifolds $L\subset M$
which are conical at infinity, i.e.\ invariant under the 
flow of the Liouville vector field outside of a compact subset, and
such that the exact 1-form $\theta_{|L}$ vanishes outside of a compact set.
In other terms, if $L$ is noncompact then at infinity it must coincide with the cone 
$(1,\infty)\times \partial L$ over some
Legendrian submanifold $\partial L$ of $\partial M$.

The Hamiltonian perturbations used to define Floer complexes in the
wrapped setting are very specific: namely, we only consider Hamiltonians
$H:M\to \R$ which, outside of a compact subset of $M$, satisfy 
$H=r^2$ where $r\in (1,\infty)$ is the radial coordinate of the
symplectization $(1,\infty)\times \partial M$. Thus, outside of a
compact set the Hamiltonian vector field $X_H$ is equal to $2r$ times the
Reeb vector field $R_\alpha$ of the contact form $\alpha$ on $\partial M$.

Given two objects $L_0,L_1$, the generating set $\mathcal{X}(L_0,L_1)$ of
the {\em wrapped Floer complex} $CW(L_0,L_1)=CW(L_0,L_1;H)$ consists of time 1 trajectories
of the flow of $X_H$ which start on $L_0$ and end on $L_1$, i.e.\ points of
$\phi^1_H(L_0)\cap L_1$. More concretely, these consist of
(perturbed) intersections between 
$L_0$ and $L_1$ in the interior $M^{in}$ on one hand, and Reeb chords
(of arbitrary length) from $\partial L_0$ to $\partial L_1$ on 
the other hand (see Figure \ref{fig:wrap}). Thus, wrapped Floer cohomology is
closely related to Legendrian contact homology.
(Of course, we need to assume that $\phi^1_H(L_0)$ intersects $L_1$
transversely, and in particular that the Reeb chords from $\partial L_0$
to $\partial L_1$ are non-degenerate; otherwise a small
modification of $H$ is required.) 

The differential on the wrapped Floer complex counts solutions to Floer's
equation \eqref{eq:floereq}, i.e.\ perturbed $J$-holomorphic strips with
boundary on $L_0$ and $L_1$, as in~\S \ref{s:HF}. (Note: due to exactness
we can work directly over the field $\mathbb{K}$, without resorting to
Novikov coefficients.) As in Remark \ref{rmk:hamtrick}, these can
equivalently
be viewed as $(\phi_H^{1-t})_*J$-holomorphic strips with boundary on $\phi_H^1(L_0)$ and
$L_1$. The assumptions made on the objects of $\mathcal{W}(M)$
and on the Hamiltonian $H$ ensure that, for suitably chosen $J$, perturbed $J$-holomorphic strips
are well-behaved: an {\em a priori} energy estimate ensures that all solutions
of \eqref{eq:floereq}
which converge to a given generator $p\in \mathcal{X}(L_0,L_1)$ as $s\to+\infty$
remain within a bounded subset of
$M$ (see e.g.\ \cite{AbGenerate}). Thus, $\partial p$ is a finite linear
combination of generators of the wrapped Floer complex.

A subtlety comes up when we attempt to define the product operation on
wrapped Floer complexes,
\begin{equation}\label{eq:wrapproduct}
CW(L_1,L_2;H)\otimes CW(L_0,L_1;H)\to CW(L_0,L_2;H).
\end{equation}
For the perturbed
Cauchy-Riemann equation \eqref{eq:perthol} to be well-behaved and satisfy
{\em a priori} energy estimates in spite of the non-compactness of $M$,
one needs the 1-form $\beta$ that appears in the perturbation term
$X_H\otimes \beta$ to satisfy
$d\beta\le 0$ (cf.\ \cite{AS,AbGenerate}). In other terms, the naturally 
defined product map would take values in $CW(L_0,L_2;2H)$, and the usual
continuation map from this complex to $CW(L_0,L_2;H)$ fails to be
well-defined. This can be remedied
using the following rescaling trick alluded to in \cite{FSS} and
systematically developed in \cite{AbGenerate}. 

Recall that 
the flow of the Liouville vector field is conformally symplectic and, 
in the symplectization $(1,\infty)\times \partial M$ where
$Z=r\,\frac{\partial}{\partial r}$, simply amounts to rescaling in the $r$ direction. 
For $\rho>1$, denote by $\psi^\rho$ the time $\log \rho$ flow of $Z$, which rescales 
$r$ by a factor of $\rho$. Then there is a natural isomorphism
\begin{equation}\label{eq:rescaletrick}
CW(L_0,L_1;H,J)\cong CW(\psi^\rho(L_0),\psi^\rho(L_1);\rho^{-1}H\circ
\psi^\rho,\psi^\rho_*J).\end{equation}
Moreover, our assumptions imply that $\psi^\rho(L_i)$ is exact Lagrangian
isotopic to $L_i$ by a compactly supported isotopy, 
and $\rho^{-1}H\circ \psi^\rho$ coincides with $\rho H$ at infinity.
Abouzaid shows that these properties ensure the existence of a well-defined product map
\begin{equation}\label{eq:productCW}
CW(L_1,L_2;H,J)\otimes CW(L_0,L_1;H,J)\to CW(\psi^2(L_0),
\psi^2(L_2);{\tfrac12}H\circ \psi^2,\psi^2_*J),\end{equation}
determined by counts of index 0 finite energy maps $u:D\to M$ from a disc with three strip-like
ends to $M$, mapping the three components of $\partial D$ to the images of the respective
Lagrangians under suitable Liouville rescalings, and solving the perturbed
Cauchy-Riemann equation
$$\bigl(du-X_{\tilde{H}}\otimes\beta\bigr)^{0,1}_{\tilde{J}}=0,$$ where $\beta$ is
a closed 1-form on $D$ with $\beta_{|\partial D}=0$ which is standard
in the strip-like ends (modelled on $dt$ for the input ends, $2\,dt$ for the
output end), and $\tilde{H}$ and $\tilde{J}$ are obtained from $H$ and $J$
by suitable rescalings ($\tilde{H}=H$ and $\tilde{J}=J$ near
the input punctures; $\tilde{H}=\frac14 H\circ \psi^2$ and $\tilde{J}=\psi^2_*J$
near the output puncture; see \cite{AbGenerate}).
The map \eqref{eq:productCW}, composed with the isomorphism \eqref{eq:rescaletrick}, 
yields the desired product map \eqref{eq:wrapproduct}. 
The higher products
$$\mu^k:CW(L_{k-1},L_k;H)\otimes \dots\otimes CW(L_0,L_1;H)\to
CW(L_0,L_k;H)$$ are constructed in the same manner \cite{AbGenerate}.
These structure maps make $\mathcal{W}(M)$ an $A_\infty$-category, the
wrapped Fukaya category of the Liouville manifold $M$. 

\begin{remark}
The rescaling trick can be informally understood as follows. As
mentioned above, the naturally defined product map on wrapped Floer
complexes takes values in $CW(L_0,L_2;2H)$; while the usual construction
of a continuation map cannot be used to map this complex to $CW(L_0,L_2;H)$,
the fact that $\frac12 H\circ \psi^2=2H$ at infinity and the assumptions
made on $L_0$ and $L_2$ imply that there is a well-defined continuation 
map to $CW(\psi^2(L_0),\psi^2(L_2);\frac12 H\circ \psi^2)$, which by \eqref{eq:rescaletrick} is isomorphic to
$CW(L_0,L_2;H)$. (Note: while this is a slightly simpler way
to describe the cohomology-level product, it lacks the compatibility and
consistency features needed to construct the chain-level
$A_\infty$-structure, hence the slightly more complicated construction in
\cite{AbGenerate}).
\end{remark}

\begin{remark}
Since compact exact Lagrangian submanifolds of $M^{in}$ are not 
affected by the wrapping at infinity, $\mathcal{W}(M)$ contains the ordinary 
Fukaya category (of compact exact Lagrangian submanifolds) 
as a full $A_\infty$-subcategory.
\end{remark}

\subsection{An example} \label{ss:wrapcyl}

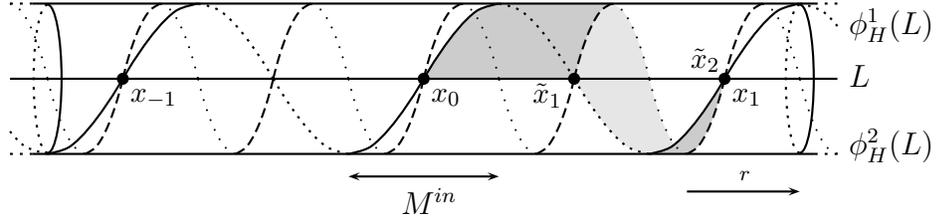
\begin{figure}[t]
\setlength{\unitlength}{1cm}
\begin{picture}(11,2.5)(-5.5,-1.5)
\psset{unit=\unitlength,dash=1pt 2.5pt}
\newgray{gray30}{0.8}
\newgray{gray15}{0.9}
\pscurve[linestyle=none,fillstyle=solid,fillcolor=gray15]%
 (3,-1)(2.7,-0.9)(2,0)(2,0)(1.3,0.9)(1,1)(1,1)(1,1)(2.5,1)(2.5,1)(2.5,1)%
 (2.65,0.9)(3,0)(3,0)(3.35,-0.9)(3.5,-1)
\pscurve[linestyle=none,fillstyle=solid,fillcolor=gray30]%
 (3,-1)(3.3,-0.9)(4,0)(4,0)(4,0)(3.65,-0.9)(3.5,-1)
\pscurve[linestyle=none,fillstyle=solid,fillcolor=gray30]%
 (0,0)(0.7,0.9)(1,1)(1,1)(1,1)(2.5,1)(2.5,1)(2.5,1)%
 (2.35,0.9)(2,0)
\psline(-5.2,-1)(5.2,-1)
\psline(-5.2,1)(5.2,1)
\psline[linestyle=dashed](-5.5,-1)(-5.2,-1)
\psline[linestyle=dashed](-5.5,1)(-5.2,1)
\psline[linestyle=dashed](5.5,-1)(5.2,-1)
\psline[linestyle=dashed](5.5,1)(5.2,1)
\psellipticarc(-5,0)(0.2,1){-90}{90}
\psellipticarc[linestyle=dashed](-5,0)(0.2,1){90}{270}
\psellipticarc(5,0)(0.2,1){-90}{90}
\psellipticarc[linestyle=dashed](5,0)(0.2,1){90}{270}
\put(-0.3,-1.8){$M^{in}$}
\psline{<->}(-1,-1.3)(1,-1.3)
\psline{->}(3.5,-1.5)(5,-1.5)
\put(4.2,-1.35){\tiny $r$}
\psline(-5.5,0)(5.5,0)
\put(5.65,-0.1){$L$}
\pscurve[linestyle=dashed](-5,-1)(-5.3,-0.9)(-5.5,-0.7)
\pscurve(-4,0)(-4.7,-0.9)(-5,-1)
\pscurve(-4,0)(-3.3,0.9)(-3,1)
\pscurve[linestyle=dashed](-2,0)(-2.7,0.9)(-3,1)
\pscurve[linestyle=dashed](-1,-1)(-1.3,-0.9)(-2,0)
\pscurve(-1,-1)(-0.7,-0.9)(0,0)
\pscurve(0,0)(0.7,0.9)(1,1)
\pscurve[linestyle=dashed](1,1)(1.3,0.9)(2,0)
\pscurve[linestyle=dashed](2,0)(2.7,-0.9)(3,-1)
\pscurve(3,-1)(3.3,-0.9)(4,0)
\pscurve(4,0)(4.7,0.9)(5,1)
\pscurve[linestyle=dashed](5,1)(5.3,0.9)(5.5,0.7)
\put(5.65,0.6){$\phi_H^1(L)$}
\pscircle*(0,0){0.08}
\pscircle*(4,0){0.08}
\pscircle*(-4,0){0.08}
\put(0.1,-0.3){\small $x_0$}
\put(4.1,-0.3){\small $x_1$}
\put(-3.9,-0.3){\small $x_{-1}$}
\pscurve[linestyle=dotted](-5,0)(-5.35,0.9)(-5.5,1)
\pscurve[linestyle=dotted](-4.5,-1)(-4.65,-0.9)(-5,0)
\pscurve[linestyle=dashed,dash=4pt 2pt](-4,0)(-4.35,-0.9)(-4.5,-1)
\pscurve[linestyle=dashed,dash=4pt 2pt](-4,0)(-3.65,0.9)(-3.5,1)
\pscurve[linestyle=dotted](-3,0)(-3.35,0.9)(-3.5,1)
\pscurve[linestyle=dotted](-2.5,-1)(-2.65,-0.9)(-3,0)
\pscurve[linestyle=dashed,dash=4pt 2pt](-2,0)(-2.35,-0.9)(-2.5,-1)
\pscurve[linestyle=dashed,dash=4pt 2pt](-2,0)(-1.65,0.9)(-1.5,1)
\pscurve[linestyle=dotted](-1,0)(-1.35,0.9)(-1.5,1)
\pscurve[linestyle=dotted](-0.5,-1)(-0.65,-0.9)(-1,0)
\pscurve[linestyle=dashed,dash=4pt 2pt](-0.5,-1)(-0.35,-0.9)(0,0)
\pscurve[linestyle=dashed,dash=4pt 2pt](0,0)(0.35,0.9)(0.5,1)
\pscurve[linestyle=dotted](0.5,1)(0.65,0.9)(1,0)
\pscurve[linestyle=dotted](1,0)(1.35,-0.9)(1.5,-1)
\pscurve[linestyle=dashed,dash=4pt 2pt](1.5,-1)(1.65,-0.9)(2,0)
\pscurve[linestyle=dashed,dash=4pt 2pt](2,0)(2.35,0.9)(2.5,1)
\pscurve[linestyle=dotted](2.5,1)(2.65,0.9)(3,0)
\pscurve[linestyle=dotted](3,0)(3.35,-0.9)(3.5,-1)
\pscurve[linestyle=dashed,dash=4pt 2pt](3.5,-1)(3.65,-0.9)(4,0)
\pscurve[linestyle=dashed,dash=4pt 2pt](4,0)(4.35,0.9)(4.5,1)
\pscurve[linestyle=dotted](4.5,1)(4.65,0.9)(5,0)
\pscurve[linestyle=dotted](5,0)(5.35,-0.9)(5.5,-1)
\put(5.65,-1){$\phi_H^2(L)$}
\pscircle*(2,0){0.08}
\put(1.45,-0.32){\small $\tilde{x}_1$}
\put(3.55,0.15){\small $\tilde{x}_2$}
\end{picture}
\caption{The wrapped Floer cohomology of $L=\R\times \{pt\}$ in $\R\times S^1$}
\label{fig:wrapcyl}
\end{figure}

Let $M=T^*S^1=\R\times S^1$, equipped with the
standard Liouville form $r\,d\theta$ and the wrapping Hamiltonian $H=r^2$, 
and consider the exact Lagrangian submanifold $L=\R\times \{pt\}$.
We can label the intersection points of $\phi_H^1(L)$ with $L$ by
integers, $\mathcal{X}(L,L)=\{x_i,\ i\in \Z\}$, in increasing order along
the real axis, where $x_0$ is the intersection occurring at the minimum of
$H$; in other terms, $x_0$ is an interior intersection of $L$ with a small
pushoff of it, while the other generators correspond to Reeb
chords from $\partial L=\{pt\}\sqcup \{pt\}$ to itself in the contact 
manifold $\partial M=S^1\sqcup S^1$ (see Figure \ref{fig:wrapcyl}).

Recall that the differential on $CW(L,L)$ counts rigid
pseudo-holomorphic strips (for a $t$-dependent almost-complex structure)
with boundary on $L$ and $\phi_H^1(L)$. Since there are no such strips
(see Figure \ref{fig:wrapcyl}), the Floer differential
on $CW(L,L)$ vanishes identically, and $HW(L,L)\simeq CW(L,L)=\mathrm{span}
\,\{x_i,\ i\in\Z\}$. (This can also be seen by observing that all generators of $CW(L,L)$ have degree $0$ for the
natural $\Z$-grading.)

The product structure on $CW(L,L)$ counts perturbed pseudo-holomorphic
discs with three strip-like ends, as explained above; in the present case,
$L$ is invariant under the Liouville flow $\psi^\rho:(r,\theta)\mapsto
(\rho r,\theta)$, while $H\circ \psi^\rho=\rho^2 H$. Thus, the
rescaling trick only affects the almost-complex structure (i.e.,
$\psi^2$ intertwines $CW(L,L;H,J)$ and $CW(L,L;2H,\psi^2_*J)$), and otherwise
simply amounts to identifying $\mathcal{X}(L,L;2H)=\phi_H^2(L)\cap L$ with
$\mathcal{X}(L,L;H)=\phi_H^1(L)\cap L$ via the radial rescaling
$r\mapsto 2r$. 

Proceeding as in Remark \ref{rmk:hamtrick}, the perturbed pseudo-holomorphic
discs with boundary on $L$ which determine the product on $CW(L,L)$ can 
then be reinterpreted as genuine pseudo-holomorphic discs (with respect
to a modified family of almost-complex structures) with boundaries on
$\phi_H^2(L)$, $\phi_H^1(L)$ and $L$. Specifically, the coefficient of a
generator $q\in \mathcal{X}(L,L)$ in the product $p_2\cdot p_1$ of two
generators $p_1,p_2\in \mathcal{X}(L,L)$ is given by a count of index 0
pseudo-holomorphic discs with boundaries on
$\phi_H^2(L)$, $\phi_H^1(L)$ and $L$, and with strip-like ends converging
to the intersection points $\phi_H^1(p_1)\in \phi_H^2(L)\cap \phi_H^1(L)$, $p_2\in \phi_H^1(L)\cap
L$, and $\tilde{q}\in \phi_H^2(L)\cap L$, where $\tilde{q}$ corresponds to
$q\in \phi_H^1(L)\cap L$ under the Liouville rescaling.

With this understood, the product structure can be determined directly by
looking at Figure \ref{fig:wrapcyl}. Observe that any two input
intersections $\phi_H^1(x_i)\in \phi_H^2(L)\cap \phi_H^1(L)$ and $x_j\in
\phi_H^1(L)\cap L$ are the vertices of a unique immersed triangle, whose
third vertex is $\tilde{x}_{i+j}\in \phi_H^2(L)\cap L$. (This is easiest
to see by lifting the diagram of Figure \ref{fig:wrapcyl} to the universal cover of
$M$.) These triangles are all regular, and we conclude that $$x_j\cdot x_i=
x_{i+j}.$$ (Recall that thanks to exactness we are working over 
$\mathbb{K}$ and not keeping track of symplectic areas.) 
For example, the triangle shaded in Figure \ref{fig:wrapcyl} illustrates the
identity $x_0\cdot x_1=x_1$. In other terms, renaming the generator $x_i$ to
$x^i$, we have a ring isomorphism 
\begin{equation}\label{eq:laurent}
CW(L,L)\simeq \mathbb{K}[x,x^{-1}].
\end{equation}
Furthermore, the higher products on $CW(L,L)$ are all identically zero, as
can be checked either by drawing the successive images of $L$ under the
wrapping flow and looking for rigid holomorphic polygons (there are none),
or more directly by recalling that $\deg(x^i)=0$ for all $i\in \Z$ whereas
$\deg(\mu^k)=2-k$. Thus \eqref{eq:laurent} is in fact an isomorphism of
$A_\infty$-algebras.

\subsection{Cotangent bundles}

The previous example is the simplest case of a general result about
cotangent bundles. Let $N$ be a compact spin manifold, and let
$M=T^*N$ equipped with its standard Liouville form $p\,dq$ and the
wrapping Hamiltonian $H=\|p\|^2$ (for some choice of Riemannian metric on
$N$). Then we have:

\begin{theorem}[Abouzaid \cite{AbCotangent1}]\label{thm:abcotangent}
Let $L=T_q^*N$, the cotangent
fiber at some point $q\in N$. Then there is a quasi-isomorphism of
$A_\infty$-algebras \begin{equation}\label{eq:abcotangent}
CW^*(L,L)\simeq C_{-*}(\Omega_qN)\end{equation}
between the wrapped Floer complex of $L=T_q^*N$ and
chains on the based loop space $\Omega_qN$ 
equipped with (an $A_\infty$-refinement of) the usual Pontryagin product.
\end{theorem}

\noindent (The corresponding statement for cohomology is an earlier result of
Abbondandolo and Schwarz \cite{AbSch}.) 

For instance, in the case of $N=S^1$, the based loop space $\Omega_q S^1$
has countably many components, each of which is contractible, thus $\Omega_q
S^1\sim \Z$, and \eqref{eq:abcotangent} reduces to~\eqref{eq:laurent}.
In fact, the assumption that $N$ is spin can be removed; in that case,
$CW^*(L,L)$ is related to chains on $\Omega_q N$ twisted by the 
$\Z$-local system determined by $w_2(N)$
\cite{AbCotangent1}.

Furthermore, Abouzaid has shown that the fiber $L=T_q^*N$ generates the
wrapped Fukaya category $\mathcal{W}(T^*N)$ \cite{AbCotangent2}. Using
Yoneda embedding (cf.\ \S \ref{ss:yoneda}), we conclude:

\begin{corollary}[Abouzaid]
The wrapped Fukaya category $\mathcal{W}(T^*N)$ quasi-embeds fully
faithfully into the category of $A_\infty$-modules over $C_{-*}(\Omega_qN)$.
\end{corollary}

\noindent (Here again, when $N$ is not spin a twist by a suitable local 
system is required.)

This and other related results can be viewed as the culmination of over a decade of
investigations of the deep connections between 
the symplectic topology of $T^*N$ and the algebraic topology of the loop
space of $N$,
as previously studied by Viterbo \cite{viterbo}, Salamon-Weber \cite{SaWe}, 
Abbondandolo-Schwarz \cite{AbSch0,AbSch}, Cieliebak-Latschev \cite{CieLa}, etc.

At the same time, studying Fukaya categories of cotangent bundles has
led to much progress on Arnold's conjecture on exact Lagrangian
submanifolds:

\begin{conj}[Arnold]\label{conj:arnexact}
Let $N$ be a compact closed manifold: then any compact closed exact Lagrangian 
submanifold of $T^*N$ (with its standard Liouville form) is Hamiltonian
isotopic to the zero section.
\end{conj}

\begin{theorem}[Fukaya-Seidel-Smith \cite{FSS}, Nadler-Zaslow \cite{NZ},
Abouzaid \cite{AbArnold}, Kragh \cite{Kragh}]
Let $L$ be a compact connected exact Lagrangian submanifold of $T^*N$. Then
as an object of $\mathcal{W}(T^*N)$, $L$ is quasi-isomorphic to the 
zero section, and the restriction of the bundle projection
$\pi_{|L}:L\to N$ is a homotopy equivalence.
\end{theorem}

Abouzaid has further shown that Floer theory detects more than purely 
topological information about exact Lagrangians in cotangent bundles:
certain exotic spheres (in dimensions $\ge 9$) do not
admit Lagrangian embeddings into $T^*S^{4k+1}$ \cite{AbSpheres}.

However, in spite of all the recent progress, Conjecture
\ref{conj:arnexact} appears to remain out of reach of current technology.

\subsection{Homological mirror symmetry}

Kontsevich's homological mirror symmetry conjecture \cite{KoICM} asserts
that the main manifestation of the phenomenon of mirror symmetry is as a
derived equivalence between the Fukaya category of a symplectic manifold and the
category of coherent sheaves of its mirror. While this conjecture
was initially stated for compact Calabi-Yau manifolds (and recently proved
for the quintic 3-fold by Sheridan \cite{Sheridan}), it also holds (and is
often easier to prove) for non-compact manifolds (in which case one should
consider the {\em wrapped} Fukaya category), and outside of the Calabi-Yau
case (in which case the mirror is a {\em Landau-Ginzburg model}, for which
one should consider Orlov's derived category of singularities
\cite{Orlov,Orlov2} rather than
the ordinary derived category of coherent sheaves).

The calculation we have performed in \S \ref{ss:wrapcyl}, together with
Abouzaid's generation statement, essentially proves homological mirror
symmetry for the cylinder $\C^*=T^*S^1$, and its mirror
$\C^*=\mathrm{Spec}\,\C[x^{\pm 1}]$. Namely, coherent sheaves over $\C^*$
are the same thing as finite rank $\C[x^{\pm 1}]$-modules. However, since
the object $L$ considered in \S \ref{ss:wrapcyl} generates the wrapped
Fukaya category, $\mathcal{W}(T^*S^1)$ quasi-embeds into the category of modules
over $CW(L,L)\simeq \C[x^{\pm 1}]$, and the image can be
characterized explicitly enough to prove the desired equivalence between
$\mathcal{W}(T^*S^1)$ and $D^b\mathrm{Coh}(\C^*)$.

This general approach extends to other examples, with the caveat that in
general there are infinitely many non-trivial higher $A_\infty$-operations;
one then needs to rely on an algebraic classification result in order to 
determine which structure coefficients need to be computed in order to
fully determine the $A_\infty$-structure up to homotopy. Symplectic
manifolds whose
Fukaya categories have
been determined in this manner include (but are not limited to)
pairs of pants \cite{AAEKO}, genus 2 curves \cite{SeGenus2}, and Calabi-Yau
hypersurfaces in projective space \cite{Sheridan}.

\subsection{An application to Heegaard-Floer homology}
Heegaard-Floer homology associates to a
closed 3-manifold $Y$ a graded abelian group $\widehat{HF}(Y)$. This
invariant is constructed
by considering a Heegaard splitting $Y=Y_1\cup_{\bar\Sigma} Y_2$ of $Y$ into two
genus $g$ handlebodies $Y_i$, each of which determines a product torus $T_i$ in the
$g$-fold symmetric product of the Heegaard surface $\bar\Sigma=\partial Y_1=
-\partial Y_2$. Deleting a marked point $z$ from $\bar\Sigma$ to obtain an open
surface $\Sigma$, $\widehat{HF}(Y)$ is then defined as the 
Floer cohomology of the Lagrangian tori $T_1,T_2$
in the symplectic manifold $\mathrm{Sym}^g(\Sigma)$, see \cite{OS}.

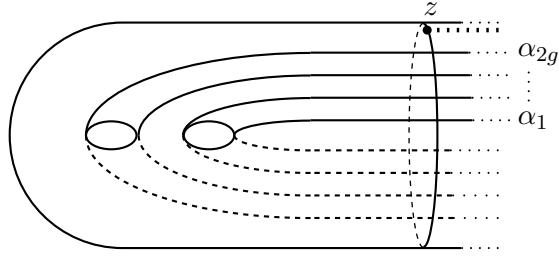
\begin{figure}[t]
\setlength{\unitlength}{1cm}
\begin{picture}(7.55,3)(-1.05,-1.5)
\psset{unit=\unitlength}
\psellipticarc[linewidth=0.5pt,linestyle=dashed,dash=2pt
2pt](4.5,0)(0.2,1.5){90}{-90}
\psellipticarc(4.5,0)(0.2,1.5){-90}{90}
\psline[linearc=1.5](5,1.5)(-1,1.5)(-1,-1.5)(5,-1.5)
\psline[linestyle=dotted](5,1.5)(5.5,1.5)
\psline[linestyle=dotted](5,-1.5)(5.5,-1.5)
\pscircle*(4.55,1.4){0.06} \put(4.5,1.6){\small $z$}
\psellipse(0.35,0)(0.35,0.2)
\psellipse(1.65,0)(0.35,0.2)
\psellipticarc(3,0)(1.03,0.21){90}{180}
\psellipticarc(3,0)(1.71,0.51){90}{180}
\psellipticarc(3,0)(2.3,0.81){90}{180}
\psellipticarc(3,0)(3,1.11){90}{180}
\psline(3,0.2)(5.15,0.2)
\psline(3,0.5)(5.14,0.5)
\psline(3,0.8)(5.12,0.8)
\psline(3,1.1)(5.1,1.1)
\psline[linestyle=dotted](5.15,0.2)(5.65,0.2)
\psline[linestyle=dotted](5.14,0.5)(5.64,0.5)
\psline[linestyle=dotted](5.12,0.8)(5.62,0.8)
\psline[linestyle=dotted](5.1,1.1)(5.6,1.1)
\psellipticarc[linestyle=dashed,dash=2pt 2pt](3,0)(1.03,0.215){180}{270}
\psellipticarc[linestyle=dashed,dash=2pt 2pt](3,0)(1.71,0.515){180}{270}
\psellipticarc[linestyle=dashed,dash=2pt 2pt](3,0)(2.3,0.815){180}{270}
\psellipticarc[linestyle=dashed,dash=2pt 2pt](3,0)(3,1.115){180}{270}
\psline[linestyle=dashed,dash=2pt 2pt](3.05,-0.2)(4.85,-0.2)
\psline[linestyle=dashed,dash=2pt 2pt](3.05,-0.5)(4.86,-0.5)
\psline[linestyle=dashed,dash=2pt 2pt](3.05,-0.8)(4.88,-0.8)
\psline[linestyle=dashed,dash=2pt 2pt](3.05,-1.1)(4.9,-1.1)
\psline[linestyle=dotted](4.85,-0.2)(5.5,-0.2)
\psline[linestyle=dotted](4.86,-0.5)(5.5,-0.5)
\psline[linestyle=dotted](4.88,-0.8)(5.5,-0.8)
\psline[linestyle=dotted](4.89,-1.1)(5.5,-1.1)
\put(5.75,0.15){\small $\alpha_1$}
\put(5.75,1.05){\small $\alpha_{2g}$}
\psline[linestyle=dotted](5.9,0.85)(5.9,0.45)
\psline[linestyle=dashed,linewidth=0.06,dash=1pt 3pt](4.55,1.4)(5.5,1.4)
\end{picture}
\caption{Generating $\mathcal{W}(\mathrm{Sym^g}(\Sigma))$}
\label{fig:symg}
\end{figure}

In this context it is natural to study the Fukaya category
(ordinary or wrapped) of $\mathrm{Sym}^g(\Sigma)$ (equipped with a K\"ahler
form which agrees with the product one away from the diagonal). It turns out that the
wrapped category has a particularly nice set of generators. Namely, consider
a collection of $2g$ disjoint properly embedded arcs
$\alpha_1,\dots,\alpha_{2g}$ in $\Sigma$ such that $\Sigma\setminus
(\alpha_1\cup\dots\cup \alpha_{2g})$ is homeomorphic to a disc, see e.g.\
Figure \ref{fig:symg}. Given a $g$-element subset $s\subseteq
\{1,\dots,2g\}$, the product $D_s=\prod_{i\in s} \alpha_i$ is an exact
Lagrangian submanifold of $\mathrm{Sym}^g(\Sigma)$, and we have:

\begin{theorem}\cite{Aur10,Aur10b}\label{thm:fuksymg}
The Lagrangian submanifolds $D_s=\prod_{i\in s} \alpha_i$, $s\subseteq
\{1,\dots,2g\}$, $|s|=g$ generate $\mathcal{W}(\mathrm{Sym}^g(\Sigma))$.
\end{theorem}

Thus, by Yoneda embedding, Lagrangian submanifolds of
$\mathrm{Sym}^g(\Sigma)$ can be viewed as modules over the
$A_\infty$-algebra $\bigoplus_{s,s'} \mathrm{hom}(D_s,D_{s'})$.

Determining this $A_\infty$-algebra is not completely hopeless, as
the wrapping Hamiltonian $H$ on $\mathrm{Sym}^g(\Sigma)$ can be chosen
in a manner compatible with the product structure so that $\phi_H^1(D_s)=
\prod_{i\in s} \phi_h^1(\alpha_i)$, where $h$ is a Hamiltonian on $\Sigma$
that grows quadratically in the cylindrical end, and pseudo-holomorphic
discs in the symmetric product can be viewed by projecting them to $\Sigma$ as is customary in
Heegaard-Floer theory; nonetheless, things are complicated by the
presence of many nontrivial $A_\infty$-products.

It is easier to study a {\em partially wrapped} version of the Fukaya
category, in which the wrapping ``stops'' along a ray $\{z\}\times
(1,\infty)$ in the cylindrical end of $\Sigma$; i.e., the
Hamiltonian is again chosen to be compatible with the product structure away from
the diagonal, but the effect on each component is to push the ends of the
arc $\alpha_i$ in the positive direction towards the ray 
$\{z\}\times (1,\infty)$, without ever crossing it: see~\cite{Aur10}.
Theorem \ref{thm:fuksymg} continues to hold in this setting: the product
Lagrangians $D_s$ also generate the partially wrapped Fukaya category.
Furthermore, in the partially wrapped case the $A_\infty$-algebra
$\mathcal{A}=\bigoplus_{s,s'} \mathrm{hom}(D_s,D_{s'})$ turns out to be a
finite-dimensional dg-algebra (i.e., $\mu^k=0$ for $k\ge 3$) which
admits a simple explicit combinatorial description \cite{Aur10}; in fact, 
$\mathcal{A}$ is precisely the
{\em strands algebra} first introduced by Lipshitz, 
Ozsv\'ath and Thurston \cite{LOT}.

By Yoneda embedding, Lagrangian submanifolds of $\mathrm{Sym}^g(\Sigma)$,
such as the product tori associated to genus $g$ handlebodies in
Heegaard-Floer theory, can be viewed as $A_\infty$-modules over the 
strands algebra. Moreover, the same holds true for {\em generalized 
Lagrangian submanifolds} of $\mathrm{Sym}^g(\Sigma)$ (i.e., 
formal images of Lagrangian submanifolds under sequences of Lagrangian
correspondences, cf.\ \cite{WW}), such as those associated to arbitrary 
3-manifolds with boundary $\bar\Sigma$ (not just handlebodies) according
to ongoing work of Lekili and Perutz. This provides a symplectic geometry
interpretation of Lipshitz-Ozsv\'ath-Thurston's {\em bordered Heegaard-Floer
homology} \cite{LOT}, which associates to a 3-manifold $Y$ with boundary $\partial Y=
\bar\Sigma$ an $A_\infty$-module $\widehat{CFA}(Y)$ over the strands
algebra. Namely, Lekili and Perutz's construction associates to such a 3-manifold
a generalized Lagrangian submanifold of $\mathrm{Sym}^g(\Sigma)$, whose
image under Yoneda embedding (as in \S \ref{ss:yoneda}, but using {\em
quilted Floer cohomology} of Lagrangian correspondences) is the
$A_\infty$-module $\widehat{CFA}(Y)$; see
\cite{Aur10,Aur10b}.

\subsection{A closing remark}
The methods available to calculate Floer cohomology and Fukaya categories are
still evolving rapidly. Besides the use of algebraic generation statements 
such as those in \cite{AbGenerate} and \cite{SeBook} to reduce to a simpler 
set of Lagrangian submanifolds, there are at least two key ideas that have 
made calculations possible.

On one hand, it is
often possible to find holomorphic projection maps (to the complex plane or
to other Riemann surfaces) under which the given Lagrangians project to arcs
or curves, in which case holomorphic discs can be studied by looking at
their projections to the base and by reducing to the symplectic geometry of
the fiber; this is e.g.\ the guiding principle of Seidel's work on
Lefschetz fibrations \cite{SeBook,SeLF} and the various calculations done
using that framework. 

At the same time, since such holomorphic projections are easier to
come by on open manifolds, another idea that nicely complements this one is
to carry out calculations for an exact open subdomain $M^0$ of the given
symplectic manifold $M$ obtained by deleting some complex hypersurface, 
and then use abstract deformation theory to view the Fukaya 
category of $M$ as an $A_\infty$-deformation of that of $M^0$ (cf.\
\cite{SeICM}). The
Hochschild cohomology class that determines the deformation is then often
determined by symmetry considerations and/or by studying specific 
$A_\infty$-structure maps (i.e., certain counts of holomorphic discs in $M$).
See e.g.\ \cite{SeGenus2,Sheridan} for an illustration of
this approach. (One guiding principle which might explain why this approach
is so successful is that algebraic deformations of Fukaya categories are
often {\em geometric}: natural ``closed-open'' maps from the 
quantum or symplectic cohomology of $M$ to the Hochschild cohomology of
its ordinary or wrapped Fukaya category often turn out to be isomorphisms
\cite{AFOOO,Ganatra}.)

Going forward, there is hope that sheaf-theoretic methods will lead to
completely new methods of computation of Fukaya categories (at least for 
Liouville manifolds) in terms of the topology of a Lagrangian ``skeleton''.
This is an idea that to our knowledge originated with Kontsevich 
\cite{KoSkel}, and was subsequently developed by
various other authors (see e.g.\ \cite{SeSkel,AbPlumbings,STZ,NadlerSkel}); the ultimate
goal being to bypass the analysis of pseudo-holomorphic curves in favor of
algebraic and topological methods. It is too early to tell how successful
these approaches will be, but it is entirely possible that
they will ultimately supplant the techniques we have described in this text.

\end{document}